\documentclass[12pt]{amsart}
\usepackage{amsmath}
\usepackage{amsthm}
\usepackage{amssymb}
\usepackage{amsfonts}
\usepackage[applemac]{inputenc}
\usepackage{mathrsfs}

\newcommand{\real}{\mathbb{R}}
\newcommand{\com}{\mathbb{C}}
\newcommand{\Z}{\mathbb{Z}}
\newcommand{\N}{\mathbb N}
\newcommand{\cfun}{\mathcal C}

\newcommand{\de}{\mathrm{d}}
\newcommand{\im}{\mathrm{Im}}
\newcommand{\re}{\mathrm{Re}}
\newcommand{\eps}{\varepsilon}
\newcommand{\wlim}{\rightharpoonup}

\newcommand{\wstar}{\overset{\ast}{\wlim}}
\newcommand{\1}[1]{\frac{1}{#1}}
\newcommand{\mez}{\frac{1}{2}}

\renewcommand{\d}{\partial}

\newcommand{\conj}[1]{\overline{#1}}
\newcommand{\tonde}[1]{\left(#1\right)}
\newcommand{\quadre}[1]{\left[#1\right]}
\newcommand{\graffe}[1]{\left\{#1\right\}}
\newcommand{\norma}[1]{ \left\|#1\right\|}

\newcommand{\slabint}{\int_{k\tau}^{(k+1)\tau}\int_{\real^3}}

\newcommand{\lplq}[4]{\|#1\|_{L^{#2}_tL^{#3}_x(#4\times\real^3)}}

\newcommand{\lplqs}[3]{\|#1\|_{L^{#2}_tL^{#3}_x}}

\newcommand{\LpLq}[4]{\norma{#1}_{L^{#2}_tL^{#3}_x\tonde{#4\times\real^3}}}

\DeclareMathOperator{\supp}{supp}
\DeclareMathOperator{\diver}{div}

\newtheorem{thm}{Theorem}
\newtheorem{prop}[thm]{Proposition}
\newtheorem{lemma}[thm]{Lemma}
\newtheorem{oss}[thm]{Remark}
\newtheorem{defn}[thm]{Definition}
\newtheorem{cor}[thm]{Corollary}

\newenvironment{dimo}
{\proof}
{\hfill $\square$\vspace{1cm}\newline}

\author{Paolo Antonelli}
\address{Paolo Antonelli --- 
Dipartimento di Matematica Pura ed Applicata \\
Universit\`a degli Studi dell'Aquila \\
Via Vetoio\\
		     67010  Coppito (AQ), Italy}
\email{paolo.antonelli@univaq.it}

\author{Pierangelo Marcati}
\address{Pierangelo Marcati  --- 
Dipartimento di Matematica Pura ed Applicata \\
Universit\`a degli Studi dell'Aquila \\
Via Vetoio\\
		     67010  Coppito (AQ), Italy}
\email{marcati@univaq.it}

\title[Finite energy weak solutions to QHD]{On the finite energy weak solutions to a system in Quantum Fluid Dynamics}
\date{\today}

\begin{document}
\begin{abstract}
In this paper we consider the global existence of weak solutions to  a
class of Quantum Hydrodynamics (QHD) systems with initial data,
arbitrarily large in the energy norm.  These type of  models,
initially proposed by Madelung \cite{M}, have been extensively used in
Physics to investigate Supefluidity and  Superconductivity 
phenomena \cite{F}, \cite{L} and more recently in the modeling of semiconductor
devices \cite{G} . Our approach is based on various tools, namely the
wave functions polar decomposition, the  construction of approximate
solution via a fractional steps method which iterates a Schr\"odinger 
Madelung  picture with a suitable wave function updating mechanism.
Therefore  several \emph{a priori} bounds of energy, dispersive and local
smoothing type allow us to prove the compactness of the approximating
sequences.  No uniqueness result is provided.
\end{abstract}
\maketitle
\section{introduction}
In this paper we study the Cauchy problem for the Quantum Hydrodynamics (QHD) system:
\bigskip
\begin{equation}\label{eq:QHD3}
\indent
\left\{\begin{array}{l}
\d_t\rho+\diver J=0\\
\d_tJ+\diver\tonde{\frac{J\otimes J}{\rho}}+\nabla P(\rho)+\rho\nabla V+f(\sqrt{\rho}, J, \nabla\sqrt{\rho})
=\frac{\hbar^2}{2}\rho\nabla\tonde{\frac{\Delta\sqrt{\rho}}{\sqrt{\rho}}}\\
-\Delta V=\rho-C(x),
\end{array}\right.
\bigskip
\end{equation}
with initial data
\begin{equation}\label{eq:QHD_IV3}
\rho(0)=\rho_0,\quad J(0)=J_0.
\end{equation}
We are interested to study the global existence in the class of finite energy initial data without higher regularity hypotheses or smallness assumptions.
\newline 
The analysis of uniqueness of weak solutions in some restricted classes will be done in a forthcoming paper.
\newline
We will discuss in particular the case $f(\sqrt{\rho}, J, \nabla\sqrt{\rho})=J$, however we are able to
treat a more general collision term, as it will be clarified in the following Remark \ref{oss:coll}. Indeed, we show that it is possible to consider a general collision term of the form $f=\alpha J+\rho\nabla g$, where $\alpha\geq0$ and 
$g$ is a nonlinear operator of $\sqrt{\rho}, J, \nabla\sqrt{\rho}$, satisfying certain Carathéodory-type conditions (see Remark \ref{oss:coll} for more precise explanations).
\newline
There is a formal analogy  between \eqref{eq:QHD3} and the classical fluid mechanics system, in particular  when $\hbar  = 0$, the system  \eqref{eq:QHD3}  formally coincides with the (nonhomogeneous) Euler-Poisson incompressible fluid system. 
\newline
The theoretical description of microphysical systems is generally based on the wave mechanics of 
Schr\"odinger , the matrix mechanics of Heisenberg or the path-integral mechanics of  Feynman (see \cite{F}).  Another approach to Quantum mechanics was taken by Madelung and  de Broglie (see \cite{M}), in particular the hydrodynamic theory of quantum mechanics  has been later extended by de Broglie (the idea of  “double solution”) and used as a scheme for quasicausal interpretation of microphysical systems. 
There is an extensive literature (see for example \cite{DGPS}, \cite{L}, \cite{K}, \cite{KB}, \cite{KD}, 
\cite{HZG}, and references therein) where  superfluidity phenomena are described by means of quantum hydrodynamic  systems. 
Furthermore, the Quantum Hydrodynamics system is well known in literature since it has been used in modeling semiconductor devices at nanometric scales (see \cite{G}). 
 The hydrodynamical formulation for quantum mechanics is quite useful with respect to other descriptions for semiconductor devices, such as those based on Wigner-Poisson or  Schr\"odinger-Poisson, since kinetic or Schr\"odinger equations are computationally very expensive. For a derivation of the QHD system we refer to \cite{AI}, \cite{DR}, \cite{JM}, \cite{JMM}, \cite{GM}, \cite{DGM1}, \cite{DGM2}.
\newline
The unknowns $\rho, J$ represent  the charge and the current densities respectively, $P(\rho)$ the classical pressure which we assume to satisfy $P(\rho)=\frac{p-1}{p+1}\rho^{\frac{p+1}{2}}$ (here and throughout the paper we assume $1\leq p\leq5$). The function $V$ is the self-consistent electric potential, given by the Poisson equation, the function $C(x)$ represents the density of the background positively charged ions. In the paper we will treat extensively the case $C(x)=0$, but all the results can be extended to more general cases. For instance we can assume $C\in W^{1, 1}(\real^3)\cap W^{1, 3}(\real^3)$. 
\newline
The term 
$\frac{\hbar^2}{2}\rho\nabla\tonde{\frac{\Delta\sqrt{\rho}}{\sqrt{\rho}}}$ can be interpreted as the quantum Bohm potential, or as a quantum correction to the pressure, indeed with some regularity assumptions we can write the dispersive term in different ways:
\begin{multline}
\frac{\hbar^2}{2}\rho\nabla\tonde{\frac{\Delta\sqrt{\rho}}{\sqrt{\rho}}}
=\frac{\hbar^2}{4}\diver(\rho\nabla^2\log\rho)\\
=\frac{\hbar^2}{4}\Delta\nabla\rho-\hbar^2\diver(\nabla\sqrt{\rho}\otimes\nabla\sqrt{\rho}).
\end{multline}
\newline
There is a formal equivalence between the system \eqref{eq:QHD3} and the following nonlinear Schr\"odinger-Poisson system:
 \begin{equation}\label{eq:NLS_arg}
 \left\{\begin{array}{l}
 i\hbar\d_t\psi+\frac{\hbar^2}{2}\Delta\psi=|\psi|^{p-1}\psi+V\psi+\tilde V\psi\\
 -\Delta V=|\psi|^2
 \end{array}\right.
 \end{equation}
 where $\tilde V=\1{2i}\log\tonde{\frac{\psi}{\conj{\psi}}}$, in particular the hydrodynamic system \eqref{eq:QHD3} can be obtained by defining $\rho=|\psi|^2$, $J=\hbar\im\tonde{\conj{\psi}\nabla\psi}$ and by computing the related balance laws.
 \newline
 This problem has to face a serious mathematical difficulty connected with the need to solve 
 \eqref{eq:NLS_arg} with the ill-posed potential $\tilde V$. Presently there are no mathematical results concerning the solutions to 
 \eqref{eq:NLS_arg}, except small perturbations around constant plane waves or local existence results under various severe restrictions (see \cite{JMR}, \cite{LL}). In the paper by Li and the second author \cite{ML} there is a global existence result for the system \eqref{eq:QHD3}, regarding small perturbations in higher Sobolev norms of subsonic stationary solutions, with periodic boundary conditions.
 \newline
 A possible way to circumvent this type of difficulty could be to develop a theory regarding wave functions taking values on Riemann manifolds but we will not pursue this direction in this paper (see also \cite{CSU}, \cite{SS}).
 \newline
 Another nontrivial problem concerning the derivation of solutions to \eqref{eq:QHD3} starting from the solutions to \eqref{eq:NLS_arg}, regards the reconstruction of the initial datum $\psi(0)$ in terms of the observables 
 $\rho(0), J(0)$. Actually this is a case of a more general important problem in physics, pointed out by Weigert in \cite{W}. He named it the \emph{Pauli problem} (since this question originated from a footnote in Pauli's article in \emph{Handbuch der Physik}, see \cite{P}), and it regards the possibility of reconstructing a pure quantum state, just by knowing a finite set of measurements of the state (in our case, the mass and current densities). Here the possible existence of nodal regions, or vacuum in fluid terms, namely where $\rho=0$, forbids in general this reconstruction in a classical way, and in any case some additional requirements (quantization rules like the Bohr-Sommerfeld rule) would be necessary. In any case, various authors (see \cite{W} and references therein) showed that knowledge of only position and momentum distribution does not specify any single state.
 \newline
 The opposite direction, namely the derivation of solutions of \eqref{eq:NLS_arg} starting from solutions of \eqref{eq:QHD3} also can face severe mathematical difficulties in various points. In particular if we prescribe $\psi(0)$, we can define $\rho(0)$ and $J(0)$, however from the evolution of the quantities $\rho(t)$ and $J(t)$, we cannot reconstruct the wave function $\psi(t)$. Furthermore, from the moment equation in \eqref{eq:QHD3} we cannot derive  the quantum eikonal equation
 \begin{equation}\label{eq:eik}
 \d_tS+\mez|\nabla S|^2+h(\rho)+V+S=\frac{\hbar^2}{2}\frac{\Delta\sqrt{\rho}}{\sqrt{\rho}}
 \end{equation}
 which is the key element to reconstruct a solution of \eqref{eq:NLS_arg} via the WKB ansatz for the wave function, $\psi=\sqrt{\rho}e^{iS/\hbar}$.
 \newline
 Similar difficulties arise when approaching with Wigner functions, which has been recently quite popular to deduce quantum fluid systems in a kinetic way (see \cite{Wig}). Even in the case we know the initial data $\rho(0), J(0)$ to be originated from a wave function $\psi(0)$, it is very difficult to show that the solutions $\rho(t), J(t)$ coincide for all times with the first and second momenta of the Wigner function obtained by solving the Wigner quantum transport equation. Moreover in our case there is also the difficulty due to the non-classical potential $\tilde V$.
 \newline
 A related question has been investigated by Bourgain, Brezis, Mironescu in several papers (see 
 \cite{B} and references therein) regarding the lifting problem for harmonic maps.
 \newline
 Another formal approach is provided by the following transport equation
 \begin{equation}
 \rho\d_tS+J\nabla S+\rho h(\rho)+\rho V+\rho S=\frac{\hbar^2}{2}\sqrt{\rho}\Delta\sqrt{\rho}
 \end{equation}
 which is obtained by multiplying equation \eqref{eq:eik} by $\rho$. 
 Unfortunately the lack of regularity of the solutions does not allow to apply even the more recent advances of the theory of transport equations \cite{DiPLio}, \cite{A} or the somehow related approach developed by Teufel and Tumulka \cite{TT} in the study of trajectories of Bohmian mechanics.
 \newline
 A related problem arises in the study of Nelson stochastic mechanics (see for instance Nelson \cite{N}, Guerra and Morato \cite{GuMo}), where the mathematical theory is based on the analysis of the velocity fields  (see Carlen \cite{Carl}). The application of this approach to our problem presents the same level of difficulty of the previously mentioned methods from the transport theory.
 \newline
 A natural framework to study the existence of the weak solutions to \eqref{eq:QHD3} is given by the space of finite energy states. Here the energy associated to the system \eqref{eq:QHD3} is given by
\begin{equation}\label{eq:en}
E(t):=\int_{\real^3}\frac{\hbar^2}{2}|\nabla\sqrt{\rho(t)}|^2+\mez|\Lambda(t)|^2+f(\rho(t))
+\mez|\nabla V(t)|^2\de x,
\end{equation}
where $\Lambda:=J/\sqrt{\rho}$, and $f(\rho)=\frac{2}{p+1}\rho^{\frac{p+1}{2}}$. The function $f(\rho)$ denotes the internal energy, which is related to the pressure through the identity 
$P(\rho)=\rho f'(\rho)-f(\rho)$.
\newline
Therefore our initial data are required to satisfy
\begin{equation}\label{eq:en_in}
E_0:=\int_{\real^3}\frac{\hbar^2}{2}|\nabla\sqrt{\rho_0}|^2+\mez|\Lambda_0|^2+f(\rho_0)
+\mez|\nabla V_0|^2\de x<\infty,
\end{equation}
or equivalently (if we have $1\leq p\leq 5$),
\begin{equation}
\sqrt{\rho_0}\in H^1(\real^3)\quad\textrm{and}\quad \Lambda_0:=J_0/\sqrt{\rho_0}\in L^2(\real^3).
\end{equation}
 \begin{defn}
 We say the pair $(\rho,J)$ is a \emph{weak solution} of the Cauchy problem \eqref{eq:QHD3}, 
 \eqref{eq:QHD_IV3} in $[0, T)\times\real^3$ with Cauchy data $(\rho_0, J_0)\in L^2(\real^3)$,
  if there exist locally integrable functions $\sqrt{\rho}, \Lambda$, such that 
$\sqrt{\rho}\in L^2_{loc}([0, T);H^1_{loc}(\real^3))$, 
 $\Lambda\in L^2_{loc}([0, T);L^2_{loc}(\real^3))$ and
  by defining $\rho:=(\sqrt{\rho})^2$, $J:=\sqrt{\rho}\Lambda$, one has
   \begin{itemize}
 \item for any test function 
  $\eta\in\cfun_0^\infty([0,T)\times\real^3)$ we have
 \begin{equation}\label{eq:QHD1}
 \int_0^T\int_{\real^3}\rho\d_t\eta+J\cdot\nabla\eta\de x\de t+\int_{\real^3}\rho_0\eta(0)\de x=0;
 \end{equation}
 \item for any test function $\zeta\in\cfun_0^\infty([0,T)\times\real^3;\real^3)$
 \begin{multline}\label{eq:QHD2}
 \int_0^T\int_{\real^3}J\cdot\d_t\zeta+\Lambda\otimes\Lambda:\nabla\zeta+P(\rho)\diver\zeta
 -\rho\nabla V\cdot\zeta-J\cdot\zeta\\
 +\hbar^2\nabla\sqrt{\rho}\otimes\nabla\sqrt{\rho}:\nabla\zeta
 -\frac{\hbar^2}{4}\rho\Delta\diver\zeta\de x\de t+\int_{\real^3}J_0\cdot\zeta(0)\de x=0;
 \end{multline}
 \item \emph{(generalized irrotationality condition)} for almost every $t\in(0, T)$, 
 \begin{equation}\label{eq:irrot}
 \nabla\wedge J=2\nabla\sqrt{\rho}\wedge\Lambda
 \end{equation}
 holds in the sense of distributions.
 \end{itemize}
 \end{defn}
 \begin{oss}
 Suppose we are in the smooth case, so that we can factorize $J=\rho u$, for some current velocity field $u$, then the last condition \eqref{eq:irrot} simply means $\rho\nabla\wedge u=0$, the current velocity $u$ is irrotational in 
 $\rho\de x$. This is why we will call it \emph{generalized irrotationality condition}.
 \end{oss}
 \begin{defn}
 We say that the weak solution $(\rho, J)$ to the Cauchy problem \eqref{eq:QHD3}, \eqref{eq:QHD_IV3} is a \emph{finite energy weak solution} (FEWS) in $[0, T)\times\real^3$, if in addition for almost every 
 $t\in[0, T)$, the energy \eqref{eq:en} is finite.
 \end{defn}
 In the sequel we restrict our attention only to FEWS with $\rho_0\in L^1(\real^3)$. The hydrodynamic structure of the system \eqref{eq:QHD3}, \eqref{eq:QHD_IV3} should not lead to conclude that the solutions behave like classical fluids. Indeed the connection with Schr\"odinger equations suggests that $\rho_0, J_0$ should in any case be seen as momenta related to some wave function $\psi_0$. The main result of this paper is the existence of FEWS by assuming the initial data 
 $\rho_0, J_0$ are momenta of some wave function $\psi_0\in H^1(\real^3)$.
 \begin{thm}[Main Theorem]\label{thm:main}
 Let $\psi_0\in H^1(\real^3)$ and let us define
 \begin{equation*}
 \rho_0:=|\psi_0|^2,\qquad J_0:=\hbar\im(\conj{\psi_0}\nabla\psi_0).
 \end{equation*}
 Then, for each $0<T<\infty$, there exists a finite energy weak solution to the QHD system \eqref{eq:QHD3} in $[0, T)\times\real^3$, 
with initial data $(\rho_0, J_0)$ defined as above.
 \end{thm}
 \begin{oss}\label{oss:coll}
As we said above, we can extend the result of the previous main Theorem \ref{thm:main} to a more general type of collision term. Indeed it is possible to consider nonlinear terms of the type 
 $f(\sqrt{\rho}, J, \nabla\sqrt{\rho})=\alpha J+\rho\nabla g(t, x, \sqrt{\rho}, \Lambda, \nabla\sqrt{\rho})$, where $\alpha\geq0$ and g satisfies the following Carathèodory type conditions
\begin{itemize}
 \item for all $(\sqrt{\rho}, \Lambda)\in[0, \infty)\times\real^3$, the function 
 $(t, x)\mapsto g(t, x, \sqrt{\rho}, \Lambda, \nabla\sqrt{\rho})$ is Lebesgue measurable;
 \item for almost all $(t, x)\in[0, \infty)\times\real^3$, the function
 $(u, v, w)\mapsto g(t, x, u, v, w)$ is continuous for  $(u, v, w)\in[0, \infty)\times\real^3\times\real^3$:
 \item there exists $C\in L^\infty([0, \infty)\times\real^3)$ such that
 \begin{equation}\label{eq:g_cond}
|g(t, x, u, v, w)|\leq C(t, x)(1+|u|^4+|v|^{4/3}+|w|^{4/3}). 
 \end{equation}
 \end{itemize}
The case $\alpha=0$ requires a slight modification of the method presented in the Section \ref{sect:fQHD}, as we are going to explain below. The case $\alpha>0$ can be treated with minor modifications of the finite difference scheme  in the Section \ref{sect:fract} defined for 
$f(\sqrt{\rho}, J, \nabla\sqrt{\rho})=J$. See the Remark \ref{oss:gen_case} below in the section 
\ref{sect:fract}.
\newline
We remark that the condition \eqref{eq:g_cond} is needed in order to ensure 
$g\in L^\infty([0, \infty); L^p(\real^3)+L^\infty(\real^3))$ as long as 
$\sqrt{\rho}\in L^\infty([0, T); H^1(\real^3))$, $\Lambda\in L^\infty([0, T); L^2(\real^3))$.
Note that with this definition of $f$ we need to change the definition given before for the weak solutions. Indeed the balance law for the current density becomes
 \begin{multline*}
 \int_0^\infty\int_{\real^3}J\cdot\d_t\zeta+\Lambda\otimes\Lambda:\nabla\zeta+P(\rho)\diver\zeta
 -\rho\nabla V\cdot\zeta-J\cdot\zeta\\
 +g(t, x, \sqrt{\rho}, \Lambda)(\nabla\rho\cdot\zeta+\rho\diver\zeta)\\
 +\hbar^2\nabla\sqrt{\rho}\otimes\nabla\sqrt{\rho}:\nabla\zeta
 -\frac{\hbar^2}{4}\rho\Delta\diver\zeta\de x\de t+\int_{\real^3}J_0\cdot\zeta(0)\de x=0.
 \end{multline*}
In the Sections \ref{sect:fract} and \ref{sect:apr} we will explain how to modify our methods to include this case (see Remarks \ref{oss:gen_case} and \ref{oss:strich_gen}).
 \end{oss}
 At the end of this Section we would like to remark that the system \eqref{eq:QHD3} is closely related to the Quantum Drift-Diffusion (QDD) equation
 \begin{equation}\label{eq:qdd}
 \d_t\rho+\diver\tonde{\rho\nabla\tonde{\frac{\hbar^2}{2}\frac{\Delta\sqrt{\rho}}{\rho}-V-h(\rho)}},
 \end{equation}
 where $V$ is again the electrostatic potential and the enthalpy $h(\rho)$ is such that 
 $\rho h'(\rho)=P'(\rho)$. Indeed if we scale the collision term $J$ in \eqref{eq:QHD3} with a relaxation time $\eps$ and we write
\begin{equation}\label{eq:QHDeps}
\left\{\begin{array}{l}
\d_t\rho^\eps+\diver J^\eps=0\\
\d_tJ^\eps+\diver\tonde{\frac{J^\eps\otimes J^\eps}{\rho^\eps}}+\nabla P(\rho^\eps)
+\rho^\eps\nabla V^\eps+\1{\eps}J^\eps
=\frac{\hbar^2}{2}\rho^\eps\nabla\tonde{\frac{\Delta\sqrt{\rho^\eps}}{\sqrt{\rho^\eps}}}\\
-\Delta V^\eps=\rho^\eps-C(x),
\end{array}\right.
\end{equation} 
then in the formal limit $\eps\to0$ we get the equation \eqref{eq:qdd}. A partial result on this limit was given by J\"ungel, Li and Matsumura in \cite{JLM}, while in \cite{GST} it is shown the global existence of non-negative variational solutions to \eqref{eq:qdd} without the electrostatic potential and the pressure term. In a forthcoming paper we deal with the relaxation limit from solutions to \eqref{eq:QHDeps} to solutions of \eqref{eq:qdd}.
 \section{Preliminaries and Notations}
\subsection{Notations}
For convenience of the reader we will set some notations which will be used in the sequel.
\newline
If $X$, $Y$ are two quantities (typically non-negative), we use $X\lesssim Y$ to denote $X\leq CY$, for some absolute constant $C>0$.
\newline
We will use the standard Lebesgue norms for complex-valued measurable functions $f:\real^d\to\com$ 
\begin{equation*}
\|f\|_{L^p(\real^d)}:=\tonde{\int_{\real^d}|f(x)|^p\de x}^{\1{p}}.
\end{equation*}
If we replace $\com$ by a Banach space $X$, we will adopt the notation
\begin{equation*}
\|f\|_{L^p(\real^d;X)}:=\tonde{\int_{\real^d}\|f(x)\|_{X}\de x}^{1/p}
\end{equation*}
to denote the norm of $f:\real^d\to X$. In particular, if X is a Lebesgue space $L^r(\real^n)$, and $d=1$, we will shorten the notation by writing
\begin{equation*}
\|f\|_{L^q_tL^r_x(I\times\real^n)}:=\tonde{\int_I\|f(t)\|_{L^r(\real^n)}^q\de t}^{1/q}
=\tonde{\int_I(\int_{\real^n}|f(t,x)|^r\de x)^{q/r}\de t}^{1/q}
\end{equation*}
to denote the mixed Lebesgue norm of $f:I\to L^r(\real^n)$; moreover, we will write 
$L^q_tL^r_x(I\times\real^n):=L^q(I;L^r(\real^n))$.
\newline
For $s\in\real$ we will define the Sobolev space $H^s(\real^n):=(1-\Delta)^{-s/2}L^2(\real^n)$.
\begin{defn}[see \cite{Caz}]
We say that $(q, r)$ is an \emph{admissible pair} of exponents if $2\leq q\leq\infty, 2\leq r\leq6$, and
\begin{equation}
\1{q}=\frac{3}{2}\tonde{\mez-\1{r}}.
\end{equation}
\end{defn}
Now we will introduce the \emph{Strichartz norms}. For more details, we refer the reader to \cite{T}, \cite{CKSTT}. Let $I\times\real^3$ be a space-time slab, we define the Strichartz norm $\dot S^0(I\times\real^3)$
\begin{equation*}
\|u\|_{\dot S^0(I\times\real^3)}:=\sup\tonde{\sum_N\|P_Nu\|_{L^q_tL^r_x(I\times\real^3)}^2}^{1/2},
\end{equation*}
where the $\sup$ is taken over all the admissible pairs $(q, r)$. Here $P_N$ denotes the Paley-Littlewood projection operator, with the sum taken over diadic numbers of the form $N=2^j$, $j\in\Z$. For any 
$k\geq1$, we can define
\begin{equation*}
\|u\|_{\dot S^k(I\times\real^3)}:=\|\nabla^ku\|_{\dot S^0(I\times\real^3)}.
\end{equation*}
Note that, from Paley-Littlewood inequality we have
\begin{equation*}
\|u\|_{L^q_tL^r_x(I\times\real^3)}\lesssim\|(\sum_N|P_Nu|^2)^{1/2}\|_{L^q_tL^r_x(I\times\real^3)}
\lesssim\tonde{\sum_N\|P_Nu\|_{L^q_tL^r_x(I\times\real^3)}^2}^{1/2},
\end{equation*}
and hence for each admissible pair of exponents, one has
\begin{equation*}
\|u\|_{L^q_tL^r_x(I\times\real^n)}\lesssim\|u\|_{\dot S^0(I\times\real^3)}.
\end{equation*}
\begin{lemma}[\cite{CKSTT}]
\begin{equation*}
\|u\|_{L^4_tL^\infty_x(I\times\real^3)}\lesssim\|u\|_{\dot S^1(I\times\real^3)}.
\end{equation*}
\end{lemma}
 \subsection{Schr\"odinger equations}
 In this paragraph we recall some important results concerning the nonlinear Schr\"odinger equations that will be used throughout the paper. First of all, we recall a global well-posedness theorem for the nonlinear Schr\"odinger-Poisson system (see \cite{Caz} and the references therein):
 \begin{thm}\label{thm:NLS}
 Let us consider the following nonlinear Schr\"odinger-Poisson system
 \begin{equation}
 \left\{\begin{array}{l}
 i\hbar\d_t\psi+\frac{\hbar^2}{2}\Delta\psi=|\psi|^{p-1}\psi+V\psi\\
 -\Delta V=|\psi|^2,
 \end{array}\right.
 \end{equation}
 with the initial datum
 \begin{equation}
 \psi(0)=\psi_0\in H^1(\real^3).
 \end{equation}
 There exists a unique globally defined strong solution 
 $\psi\in\cfun(\real;H^1(\real^3))$, which depends continuously on the initial data. Furthermore, the total energy
 \begin{equation}
 E(t):=\int_{\real^3}\frac{\hbar^2}{2}|\nabla\psi(t,x)|^2+\frac{2}{p+1}|\psi(t,x)|^{p+1}
 +V(t,x)|\psi(t,x)|^2\de x
 \end{equation}
 is conserved, namely
 \begin{equation}
 E(t)=E_0,\qquad\textrm{for all }\; t\in\real.
 \end{equation}
 \end{thm}
 Moreover, the dispersive nature of Schr\"odinger equation provides also some further integrability and regularity properties of the solutions. The first result in this direction regards some space-time integrability properties,
the most important of them are  the well known \emph{Strichartz estimates} for the Schr\"odinger equation in $\real^{1+3}$. We refer to the classical paper of Ginibre and Velo \cite{GV}, the paper of Keel and Tao \cite{KT}, the books of Cazenave \cite{Caz} and Tao \cite{T} and the references therein. Let 
$\psi$ is the unique (strong) solution of the Cauchy problem for the free Schr\"odinger equation
  \begin{equation*}
  \left\{\begin{array}{l}
  i\d_t\psi+\Delta\psi=0\\
  \psi(0)=\psi_0,
  \end{array}\right.
  \end{equation*}
  in the following we shalll denote by $U(\cdot)$ the free Schr\"odinger group, defined by the indentity 
  $U(t)\psi_0=\psi(t)$.
  \newline
  The next result is taken from \cite{KT} and the estimates are obtained in a very general setting.
 \begin{thm}[Keel, Tao \cite{KT}]\label{thm:strich}
 Let $(q,r), (\tilde q, \tilde r)$ be two arbitrary admissible pairs of exponents, and let $U(\cdot)$ be the free Schr\"odinger group. Then we have
 \begin{align}
 \lplqs{U(t)f}{q}{r}&\lesssim\|f\|_{L^2(\real^3)}\\
 \lplqs{\int_{s<t}U(t-s)F(s)\de s}{q}{r}&\lesssim\lplqs{F}{\tilde q'}{\tilde r'}\\
 \lplqs{\int U(t-s)F(s)\de s}{q}{r}&\lesssim\lplqs{F}{\tilde q'}{\tilde r'}.
 \end{align}
 \end{thm}
By using the previous Theorem, we can state the following Lemma about some further integrability properties for solutions of the Schr\"odinger equations.
 \begin{lemma}[\cite{CKSTT}]
Let $I$ be a compact interval, and let $u:I\times\real^3\to\com$ be a Schwartz solution to the 
Schr\"odinger equation
\begin{equation*}
i\d_tu+\Delta u=F_1+\dotsc+F_M,
\end{equation*}
for some Schwartz functions $F_1, \dotsc, F_M$. Then we have
\begin{equation*}
\|u\|_{\dot S^0(I\times\real^n)}\lesssim\|u(t_0)\|_{L^2(\real^n)}+\lplq{F_1}{q_1'}{r_1'}{I}+\dotsc
+\lplq{F_M}{q_M'}{r_M'}{I},
\end{equation*}
where $(q_1, r_1),\dotsc, (q_M, r_M)$ are arbitrary admissible pairs of exponents.
\end{lemma}
 Furthermore, the solutions to the Schr\"odinger equation enjoy some nice local smoothing properties; there are various results in the literature regarding this property, here we recall a theorem, due to Constatin and Saut \cite{CS} which actually covers a more general setting. In particular it can been shown the solutions of the Schr\"odinger equation are locally more regular than their initial data. Close results have been proved by P. Sj\"olin \cite{Sj} and L. Vega \cite{V}.
 \begin{thm}[Constantin, Saut \cite{CS}]\label{thm:smooth1}
Let $u$ solves the free Schr\"odinger equation
\begin{equation}
\left\{\begin{array}{l}
i\d_tu+\Delta u=0\\
u(0)=u_0.
\end{array}\right.
\end{equation}
Let $\chi\in\cfun^\infty_0(\real^{1+3})$ of the form
\begin{equation*}
\chi(t,x)=\chi_0(t)\chi_1(x_1)\chi_2(x_2)\chi(x_3),
\end{equation*}
with $\chi_j\in\cfun^\infty_0(\real)$. Then we have
\begin{equation*}
\int_{\real^{1+3}}\chi^2(t,x)|(I-\Delta)^{1/4}u(t,x)|^2\de x\de t\leq C^2\|u_0\|_{L^2{\real^3}}.
\end{equation*}
In particular, if $u_0\in L^2(\real^3)$, one has for all $T>0$ 
\begin{equation*}
u\in L^2([0,T];H^{1/2}_{loc}(\real^3)).
\end{equation*}
\end{thm}
We have a similar result also for the nonhomogeneous case:
\begin{thm}[Constantin, Saut \cite{CS}]\label{thm:smooth2}
Let $u$ be the solution of 
\begin{equation}
\left\{\begin{array}{l}
i\d_tu+\Delta u=F\\
u(0)=u_0\in L^2(\real^3),
\end{array}\right.
\end{equation}
where $F\in L^1([0,T];L^2(\real^3))$. Then it follows
\begin{equation*}
u\in L^2([0,T];H^{1/2}_{loc}(\real^3)).
\end{equation*}
Moreover, let $\chi\in\cfun_0^\infty(\real^{1+3})$ be of the form
\begin{equation*}
\chi(t, x)=\chi_0(t)\chi_1(x_1)\chi_2(x_2)\chi_3(x_3)
\end{equation*}
with $\chi_j\in\cfun^\infty_0(\real)$, $\supp\chi_0\subset[0, T]$. Then the following local smoothing estimate holds
\begin{multline}
\tonde{\int_{\real^{1+3}}\chi^2(t,x)|(I-\Delta)^{1/4}u(t,x)|^2\de x\de t}^{\mez}\\
\leq C\tonde{\|u_0\|_{L^2(\real^3)}+\lplq{F}{1}{2}{[0,T]}}
\end{multline}
\end{thm}
This results imply that, the free Schr\"odinger group $U(\cdot)$ fulfills the following inequalities
\begin{gather}
\|U(\cdot)u_0\|_{L^2([0,T];H^{1/2}_{loc})}\lesssim\|u_0\|_{L^2(\real^3)}\label{eq:smooth1}\\
\|\int_0^tU(t-s)F(s)\de s\|_{L^2([0,T];H^{1/2}_{loc}(\real^3))}\lesssim\|F\|_{L^1([0,T];L^2(\real^3))}.
\label{eq:smooth2}
\end{gather}
\subsection{Compactness tools}
Here we recall some compactness theorems in function spaces, which will be relevant in the Section 
\ref{sect:apr} to prove the convergence of the approximate solutions. Let us recall in particular a  compactness result due to Rakotoson, Temam \cite{RT} in the spirit of classical results of Aubin, Lions and Simon, see \cite{S}, \cite{Aub}, \cite{Lions}.
\begin{thm}[Rakotoson, Temam \cite{RT}]\label{thm:comp}
Let $(V,\|\cdot\|_{V})$, $(H;\|\cdot\|_{H})$ be two separable Hilbert spaces. Assume that $V\subset H$ with a compact and dense embedding. Consider a sequence $\{u^\eps\}$, converging weakly to a function $u$ in $L^2([0,T];V)$, $T<\infty$. Then $u^\eps$ converges strongly to $u$ in 
$L^2([0,T];H)$, if and only if
\begin{enumerate}
\item $u^\eps(t)$ converges to $u(t)$ weakly in $H$ for a.e. $t$;
\item $\lim_{|E|\to0,E\subset[0,T]}\sup_{\eps>0}\int_E\|u^\eps(t)\|_H^2\de t=0$.
\end{enumerate}
\end{thm}
 \section{Polar Decomposition}\label{sect:polar}
 In this section we will explain how to decompose an arbitrary wave function $\psi$ into its amplitude 
 $\sqrt{\rho}=|\psi|$ and its unitary factor $\phi$, namely a function taking its values in the unitary circle of the complex plane, such that 
 $\psi=\sqrt{\rho}\phi$. The idea is similar in the spirit, to that one used by Y. Brenier in \cite{Br}, to find the measure preserving maps needed to write vector-valued functions as compositions of gradients of convex functions and measure preserving maps. Our case is much simpler than \cite{Br} and it can be studied directly in a simpler setting. Brenier's idea looks for projections of $L^2$ functions $u$ onto the set of measure preserving maps $S$, contained in a given sphere of $L^2$, i.e. one has to find $s\in S$, which minimizes the distance $\|u-s\|_{L^2}$, or equivalently, which maximizes $(u,s)_{L^2}$. In our case we should maximize $\re(u, s)_{L^2}$, within complex-valued functions with the constraint to take value in the unit ball of $L^\infty$. Actually in this case the maximizer turns out to be trivially determined and the previous variational argument has only to be considered a motivation to the subsequent considerations.
 \newline
 However in the case we assume $u\in H^1(\real^3)$ (that is the relevant situation in physics), then it would be useful to know the regularity of the maximizers in the light of a possible use of the considerations explained in the Remark \ref{rmk:25}.
 \newline
 Let us consider a wave function $\psi\in L^2(\real^3)$ and define the set
 \begin{equation}
 P(\psi):=\graffe{\phi\;\textrm{measurable}\;|\;\|\phi\|_{L^\infty(\real^3)}\leq1,\;\sqrt{\rho}\phi=\psi\;
 \textrm{a. e. in }\;\real^3},
 \end{equation}
 where $\sqrt{\rho}=|\psi|$. Of course if we consider $\phi\in P(\psi)$, then by the definition of the set 
 $P(\psi)$ it is immediate that $|\phi|=1$ a.e. $\sqrt{\rho}-\de x$ in $\real^3$ and $\phi\in P(\psi)$ is uniquely determined a.e. $\sqrt{\rho}-\de x$ in $\real^3$.
 \begin{oss}
 Let $B_R$ be the ball in $\real^3$ with radius $R$, centered at the origin and let us define the set
 \begin{equation*}
 P_R(\psi):=\graffe{\phi\;\textrm{measurable}\;|\:\|\phi\|_{L^\infty(\real^3)}\leq1,\;\sqrt{\rho}\phi=\psi\;
 \textrm{a. e. in }\;\real^3}.
 \end{equation*}
 Then it is easy to see that $\phi\in P_R(\psi)$ if and only if $\phi$ is a maximizer for the functional
 \begin{equation*}
 \Phi_R[\phi]:=\re\int_{B_R}\conj{\psi}\phi\de x
 \end{equation*}
 over the set
 \begin{equation*}
S_R:=\graffe{\phi\in L^2(B_R)\;|\;\|\phi\|_{L^\infty(B_R)}\leq1}.
\end{equation*}
 \end{oss}
The next lemma connects the structure of the bilinear term 
$\hbar^2\re(\nabla\conj{\psi\otimes\nabla\psi})$ with $\nabla\sqrt{\rho}\otimes\nabla\sqrt{\rho}$ and 
$\frac{J\otimes J}{\rho}$. Moreover it shows this structure is $H^1$-stable.
\begin{lemma}\label{lemma:phase}
Let $\psi\in H^1(\real^3)$, $\sqrt{\rho}:=|\psi|$, then there exists $\phi\in L^\infty(\real^3)$ such that 
$\psi=\sqrt{\rho}\phi$ a.e. in $\real^3$, $\sqrt{\rho}\in H^1(\real^3)$, 
$\nabla\sqrt{\rho}=\re(\conj{\phi}\nabla\psi)$. If we set $\Lambda:=\hbar\im(\conj{\phi}\nabla\psi)$,  one has $\Lambda\in L^2(\real^3)$ and moreover the following identity holds
\begin{equation}\label{eq:null-form}
\hbar^2\re(\d_j\conj{\psi}\d_k\psi)=\hbar^2\d_j\sqrt{\rho}\d_k\sqrt{\rho}+\Lambda^{(j)}\Lambda^{(k)}.
\end{equation}
Furthermore, let $\psi_n\to\psi$ strongly in $H^1(\real^3)$, then  it follows
\begin{equation}\label{eq:strong_conv}
\nabla\sqrt{\rho_n}\to\nabla\sqrt{\rho},\quad \Lambda_n\to\Lambda\qquad\textrm{in}\; L^2(\real^3),
\end{equation}
where $\Lambda_n:=\hbar\im(\conj{\phi_n}\nabla\psi_n)$.
\end{lemma}
\begin{dimo}
Let us consider a sequence $\{\psi_n\}\subset\cfun^\infty(\real^3)$, $\psi\to\psi$ in 
$H^1(\real^3)$, define as before
\begin{equation}\label{defn:reg-phas}
\phi_n(x)=\left\{\begin{array}{rl}
\frac{\psi_n(x)}{|\psi_n(x)|}&\textrm{if}\;\psi_n(x)\neq0\\
0&\textrm{if}\;\psi(x)=0.
\end{array}\right.
\end{equation}
Then, there exists $\phi\in L^\infty(\real^3)$ such that $\phi_n\wstar\phi$ in $L^\infty(\real^3)$ and 
$\nabla\psi_n\to\nabla\psi$ in $L^2(\real^3)$, hence
\begin{equation*}
\re\tonde{\conj{\phi_n}\nabla\psi_n}\wlim\re\tonde{\conj{\phi}\nabla\psi}\qquad\textrm{in}\;L^2(\real^3).
\end{equation*}
Since by \eqref{defn:reg-phas}, one has
\begin{equation*}
\re\tonde{\conj{\phi_n}\nabla\psi_n}=\nabla|\psi_n|
\end{equation*}
a.e. in $\real^3$, it follows
\begin{equation*}
\nabla\sqrt{\rho_n}\wlim\re\tonde{\conj{\phi}\nabla\psi}\qquad\textrm{in}\;L^2(\real^3).
\end{equation*}
Moreover, one has $\nabla\sqrt{\rho_n}\wlim\nabla\sqrt{\rho}$ in $L^2(\real^3)$, therefore
\begin{equation*}
\nabla\sqrt{\rho}=\re\tonde{\conj{\phi}\nabla\psi},
\end{equation*}
where $\phi$ is a unitary factor of $\psi$.
\newline
The identity \eqref{eq:null-form} follows immediately fro the following
\begin{multline*}
\hbar^2\re(\d_j\conj{\psi}\d_k\psi)=\hbar^2\re((\phi\d_j\conj{\psi})(\conj{\phi}\d_k\psi))
=\hbar^2\re(\phi\d_j\conj{\phi})\re(\conj{\phi}\d_k\psi)\\
-\hbar^2\im(\phi\d_j\conj{\psi})\im(\conj{\phi}\d_k\psi)
=\hbar^2\d_j\sqrt{\rho}\d_k\sqrt{\rho}+\Lambda^{(j)}\Lambda^{(k)}.
\end{multline*}
Now we prove \eqref{eq:strong_conv}. Let us take a sequence $\psi_n\to\psi$ strongly in $H^1(\real^3)$, and consider $\nabla\sqrt{\rho_n}=\re\tonde{\conj{\phi_n}\nabla\psi_n}$, 
$\Lambda_n:=\hbar\im\tonde{\conj{\phi_n}\nabla\psi_n}$. 
As before, $\phi_n\wstar\phi$ in $L^\infty(\real^3)$, with $\phi$ a polar factor for $\psi$; then $\nabla\sqrt{\rho_n}\wlim\nabla\sqrt{\rho}$, 
$\re\tonde{\conj{\phi_n}\nabla\psi_n}\wlim\re\tonde{\conj{\phi}\nabla\psi}$, and 
$\nabla\sqrt{\rho}=\re\tonde{\conj{\phi}\nabla\psi}$. Moreover, 
$\Lambda_n:=\hbar\im\tonde{\conj{\phi_n}\nabla\psi_n}\wlim\hbar\im\tonde{\conj{\phi}\nabla\psi}
=:\Lambda$.
To upgrade the weak convergence into the strong one, simply notice that by \eqref{eq:null-form} one has
\begin{multline*}
\hbar^2\|\nabla\psi\|_{L^2}^2=\hbar^2\|\nabla\sqrt{\rho}\|_{L^2}^2+\|\Lambda\|_{L^2}^2
\leq\liminf_{n\to\infty}\tonde{\hbar^2\|\nabla\sqrt{\rho_n}\|_{L^2}^2+\|\Lambda_n\|_{L^2}^2}\\
=\liminf_{n\to\infty}\hbar^2\|\nabla\psi_n\|_{L^2}^2=\hbar^2\|\nabla\psi\|_{L^2}^2.
\end{multline*}
\end{dimo}
\begin{cor}\label{cor:irrot}
Let $\psi\in H^1(\real^3)$, then
\begin{equation}\label{eq:119}
\hbar\nabla\conj{\psi}\wedge\nabla\psi=2i\nabla\sqrt{\rho}\wedge\Lambda.
\end{equation}
\end{cor}
\begin{dimo}
It suffices to note that we can write 
\begin{equation*}
\hbar\nabla\conj{\psi}\wedge\nabla\psi=\hbar(\phi\nabla\conj{\psi})\wedge(\conj{\phi}\nabla\psi),
\end{equation*}
where $\phi\in L^\infty(\real^3)$ is the polar factor of $\psi$ such that 
$\nabla\sqrt{\rho}=\re(\conj{\phi}\nabla\psi)$, $\Lambda:=\hbar\im(\conj{\phi}\nabla\psi)$, as in the Lemma 
\ref{lemma:phase}. By splitting 
$\conj{\phi}\nabla\psi$ into its real and imaginary part, we get the identity \eqref{eq:119}.
\end{dimo}
Now we state a technical lemma which will be used in the next sections. It summarizes the results of this section we use later on and will be handy for the application to the fractional step method.
\begin{lemma}\label{lemma:updat}
Let $\psi\in H^1(\real^3)$, and let $\tau, \eps>0$ be two arbitrary (small) real numbers. Then there exists 
$\tilde\psi\in H^1(\real^3)$ such that
\begin{align*}
\tilde\rho&=\rho\\
\tilde\Lambda&=(1-\tau)\Lambda+r_\eps,
\end{align*}
where $\sqrt{\rho}:=|\psi|, \sqrt{\tilde\rho}:=|\tilde\psi|, \Lambda:=\hbar\im(\conj{\phi}\nabla\psi), 
\tilde\Lambda:=\hbar\im(\conj{\tilde\phi}\nabla\tilde\psi)$, $\phi, \tilde\phi$ are polar factors for $\psi, \tilde\psi$, respectively, and
\begin{equation*}
\|r_\eps\|_{L^2(\real^3)}\leq\eps.
\end{equation*}
Furthermore we have
\begin{equation}\label{eq:121}
\nabla\tilde\psi=\nabla\psi-i\frac{\tau}{\hbar}\phi^\star\Lambda+r_{\eps, \tau},
\end{equation}
where $\|\phi^\star\|_{L^\infty(\real^3)}\leq1$ and
\begin{equation}\label{eq:122}
\|r_{\eps, \tau}\|_{L^2(\real^3)}\leq C(\tau\|\nabla\psi\|_{L^2(\real^3)}+\eps).
\end{equation}
\end{lemma}
\begin{dimo}
Consider a sequence $\{\psi_n\}\subset\cfun^\infty(\real^3)$ converging to $\psi$ in $H^1(\real^3)$, and define
\begin{equation*}
\phi_n(x):=\left\{\begin{array}{ll}
\frac{\psi_n(x)}{|\psi_n(x)|}&\textrm{if}\;\psi_n(x)\neq0\\
0&\textrm{if}\;\psi_n(x)=0
\end{array}\right.
\end{equation*}
as polar factors for the wave functions $\psi_n$. Since $\psi_n\in\cfun^\infty$, then $\phi_n$ is piecewise smooth, and $\Omega_n:=\{x\in\real^3\;:\;|\psi_n(x)|>0\}$  is an open set, with smooth boundary. Therefore we can say there exists a function $\theta_n:\Omega_n\to[0,2\pi)$, piecewise smooth in 
$\Omega_n$ and
\begin{equation*}
\phi_n(x)=e^{i\theta_n(x)},\qquad\textrm{for any}\;x\in\Omega_n.
\end{equation*}
Moreover, by the previous Lemma, we have $\phi_n\wstar\phi$ in $L^\infty(\real^3)$, where $\phi$ is a polar factor of $\psi$, and 
$\Lambda_n:=\hbar\im(\conj{\phi_n}\nabla\psi_n)\to\Lambda:=\hbar\im(\conj{\phi}\nabla\psi)$ in 
$L^2(\real^3)$. Thus there exists $n\in\N$ such that
\begin{equation*}
\|\psi-\psi_n\|_{H^1(\real^3)}+\|\Lambda-\Lambda_n\|_{L^2(\real^3)}\leq\eps.
\end{equation*}
Now we can define
\begin{equation}\label{eq:123}
\tilde\psi:=e^{i(1-\tau)\theta_n}\sqrt{\rho_n}.
\end{equation}
Let us remark that \eqref{eq:123} is a good definition for $\tilde\psi$, since outside $\Omega_n$, where $\theta_n$ is not defined, we have $\sqrt{\rho_n}=0$, and hence also $\psi_n=0$. Furthermore, 
\begin{align*}
\nabla\tilde\psi
=&e^{-i\tau\theta_n}\nabla\psi_n-i\frac{\tau}{\hbar}e^{i(1-\tau)\theta_n}\Lambda_n\\
=&\nabla\psi-i\frac{\tau}{\hbar}e^{i(1-\tau)\theta_n}\Lambda_n
+\tau\tonde{\sum_{j=1}^\infty\frac{(-i\theta_n)^j\tau^{j-1}}{j!}}\nabla\psi_n\\
&+(\nabla\psi_n-\nabla\psi)-i\frac{\tau}{\hbar}e^{i(1-\tau)\theta_n}(\Lambda_n-\Lambda)
\end{align*}
and thus we obtain \eqref{eq:121} and \eqref{eq:122}, with $r_{\eps, \tau}$ given by
\begin{equation*}
r_{\eps, \tau}=\tau\tonde{\sum_{j=1}^\infty\frac{(-i\theta_n)^j\tau^{j-1}}{j!}}\nabla\psi_n
+(\nabla\psi_n-\nabla\psi)-i\frac{\tau}{\hbar}e^{i(1-\tau)\theta_n}(\Lambda_n-\Lambda)
\end{equation*} 
Moreover
\begin{equation}
\tilde\Lambda=\hbar\im(e^{-i(1-\tau)\theta_n}\nabla\tilde\psi)
=(1-\tau)\Lambda_n
=(1-\tau)\Lambda+(1-\tau)(\Lambda_n-\Lambda)
\end{equation}
and clearly $r_\eps:=(1-\tau)(\Lambda_n-\Lambda)$ has $L^2$-norm less than $\eps$ by assumption.
\end{dimo}
\section{QHD without collisions}\label{sect:fQHD}
Now let us summarize some key points regarding the existence of weak solutions to the Quantum Hydrodynamic system, in the collisionless case.
\newline
The balance equations can be written in the following way:
\begin{equation}\label{eq:QHD_farl}
\left\{\begin{array}{ll}
\d_t\rho+\diver J=0\\
\d_tJ+\diver\tonde{\frac{J\otimes J}{\rho}}+\nabla P(\rho)+\rho\nabla V
=\frac{\hbar^2}{2}\rho\nabla\tonde{\frac{\Delta\sqrt{\rho}}{\sqrt{\rho}}}\\
-\Delta V=\rho,
\end{array}\right.
\end{equation}
where $P(\rho)=\frac{p-1}{p+1}\rho^{(p+1)/2}$, $1\leq p<5$.
\begin{defn}
We say the pair $(\rho, J)$ is a \emph{weak solution} in $[0, T)\times\real^3$ to the system 
\eqref{eq:QHD_farl} with Cauchy data $(\rho_0, J_0)\in L^1_{loc}(\real^3)$ if and only if there exist locally integrable functions 
$\sqrt{\rho}, \Lambda$ such that
$\sqrt{\rho}\in L^2_{loc}([0, T);H^1_{loc}(\real^3)), \Lambda\in L^2_{loc}([0, T);L^2_{loc}(\real^3))$ and
if we define $\rho:=(\sqrt{\rho})^2$, $J:=\sqrt{\rho}\Lambda$, then
\begin{itemize}
\item for any test function $\varphi\in\cfun^\infty_0([0, T)\times\real^3)$ we have
\begin{equation}\label{eq:w-QHD1}
\int_0^T\int_{\real^3}\rho\d_t\varphi+J\cdot\nabla\varphi\de x\de t+\int_{\real^3}\rho_0\varphi(0)\de x=0
\end{equation}
and for any test function $\eta\in\cfun^\infty_0([0, T)\times\real^3;\real^3)$ we have
\begin{multline}\label{eq:w-QHD2}
\int_0^T\int_{\real^3}J\cdot\d_t\eta+\Lambda\otimes\Lambda:\nabla\eta+P(\rho)\diver\eta
-\rho\nabla V\cdot\eta\\
+\hbar^2\nabla\sqrt{\rho}\otimes\nabla\sqrt{\rho}:\nabla\eta
-\frac{\hbar^2}{4}\rho\Delta\diver\eta\de x\de t+\int_{\real^3}J_0\cdot\eta(0)\de x=0;
\end{multline}
\item \emph{(generalized irrotationality condition)} for almost every $t\in(0, T)$
\begin{equation*}
\nabla\wedge J=2\nabla\sqrt{\rho}\wedge\Lambda
\end{equation*}
holds in the sense of distributions.
\end{itemize}
We say that $(\rho, J)$ is a \emph{finite energy weak solution} if and only if it is a weak solution and  the energy $E(t)$ is finite a.e. in $[0, T)$.
\end{defn}	
The next existence result roughly speaking shows how to get a weak solution to the system 
\eqref{eq:QHD_farl} out of a strong solution to the Schr\"odinger-Poisson system
\begin{equation}\label{eq:NLS}
\left\{\begin{array}{l}
i\hbar\d_t\psi+\frac{\hbar^2}{2}\Delta\psi=|\psi|^{p-1}\psi+V\psi\\
-\Delta V=|\psi|^2
\end{array}\right.
\end{equation}
The quadratic nonlinearities in \eqref{eq:QHD_farl} are originated by a term of the form 
$\re\tonde{\nabla\conj{\psi}\otimes\nabla\psi}$ since formally 
\begin{equation*}
\hbar^2\re\tonde{\nabla\conj{\psi}\otimes\nabla\psi}=\hbar^2\nabla\sqrt{\rho}\otimes\nabla\sqrt{\rho}
+\frac{J\otimes J}{\rho}.
\end{equation*}
However this identity can be justified in the nodal region $\{\rho=0\}$ only by means of the polar factorization discussed in the previous section.
\newline
Indeed, we stress that in the previous identity the right hand side is written in terms of $\rho$ and $J$ which do exist in the whole space $\real^3$, via the Madelung transform. However, the term 
$\frac{J\otimes J}{\rho}$ should be interpreted as $\Lambda\otimes\Lambda$, where $\Lambda$ is the Radon-Nykodim derivative of $J\de x$ with respect to $\sqrt{\rho}\de x$. Unfortunately the Madelung transformations are unable to define $\Lambda$ on the whole $\real^3$, hence one needs to use the polar factorization discussed in the previous section to define $\Lambda$ in the whole $\real^3$.
\newline
Furthermore, the  study of the existence of weak solutions of \eqref{eq:QHD_farl} is done with Cauchy data of the form $(\rho_0,J_0)=(|\psi_0|^2,\hbar\im(\conj{\psi_0}\nabla\psi_0))$, for some 
$\psi_0\in H^1(\real^3;\com)$. 
These special initial data yield to consider the Cauchy problem for the QHD system only for solutions compatible with a \emph{wave mechanics} point of view.
\begin{prop}\label{prop:19}
Let $0<T<\infty$, let $\psi_0\in H^1(\real^3)$ and define the initial data for \eqref{eq:QHD_farl}, 
$(\rho_0, J_0):=(|\psi_0|^2, \hbar\im(\conj{\psi_0}\nabla\psi_0))$. Then there exists a finite energy weak solution $(\rho, J)$ to the Cauchy problem \eqref{eq:QHD_farl} in the space-time slab $[0, T)\times\real^3$. Furthermore the energy $E(t)$ defined in \eqref{eq:en} is conserved for all times $t\in[0, T)$.
\end{prop}
The idea behind the proof of this Proposition is the following. 
Let us consider the Schr\"odinger-Poisson system \eqref{eq:NLS} with initial datum $\psi(0)=\psi_0$.
It is well known that it is globally well-posed for initial data in $H^1(\real^3)$ (see \cite{Caz}), and the solution satisfies 
$\psi\in\cfun^0(\real;H^1(\real^3))$.
\newline
Then it makes sense to define for each time $t\in[0, T)$ the quantities 
$\rho(t):=|\psi(t)|^2, J(t):=\hbar\im(\conj{\psi(t)}\nabla\psi(t))$ and we can see that $(\rho, J)$ is a finite energy weak solution of \eqref{eq:QHD_farl}: indeed for the solution $\psi$ of \eqref{eq:NLS} we have
\begin{equation*}
\d_t\nabla\psi=\frac{i\hbar}{2}\Delta\nabla\psi-\frac{i}{\hbar}\nabla\tonde{(|\psi|^{p-1}+V)\psi}
\end{equation*}
in the sense of distributions. Thus formally we get the following identities for $(\rho, J)$
\begin{equation*}
\d_t\rho=-\hbar\im(\conj{\psi}\Delta\psi)=-\diver\tonde{\hbar\im(\conj{\psi}\nabla\psi)},
\end{equation*}
\begin{multline*}
\d_tJ
=\hbar\im\tonde{\nabla\psi\tonde{-\frac{i\hbar}{2}\Delta\conj{\psi}+\frac{i}{\hbar}(|\psi|^{p-1}
+V)\conj{\psi}}}\\
+\hbar\im\tonde{\conj{\psi}\tonde{\frac{i\hbar}{2}\Delta\nabla\psi
-\frac{i}{\hbar}\nabla(|\psi|^{p-1}+V)\psi
-\frac{i}{\hbar}(|\psi|^{p-1}+V)\nabla\psi}}\\
=\frac{\hbar^2}{4}\nabla\Delta|\psi|^2
-\hbar^2\diver\tonde{\re(\nabla\conj{\psi}\otimes\nabla\psi)}
-\frac{p-1}{p+1}\nabla(|\psi|^{p+1})-|\psi|^2\nabla V.
\end{multline*}
Thanks to the Lemma \ref{lemma:phase} we can write 
\begin{equation*}
\hbar^2\re\tonde{\nabla\conj{\psi}\otimes\nabla\psi}=\hbar^2\nabla\sqrt{\rho}\otimes\nabla\sqrt{\rho}
+\Lambda\otimes\Lambda
\end{equation*}
with $\sqrt{\rho}, \Lambda$ as in statement of the Lemma \ref{lemma:phase}, and clearly they also satisfy
\begin{equation*}
J=\sqrt{\rho}\Lambda.
\end{equation*}
Hence formally the following identity holds:
\begin{equation}\label{eq:curr_dens}
\d_tJ+\diver(\Lambda\otimes\Lambda)+\nabla P(\rho)+\rho\nabla V
=\frac{\hbar^2}{4}\Delta\nabla\rho-\hbar^2\diver(\nabla\sqrt{\rho}\otimes\nabla\sqrt{\rho}).
\end{equation}
Of course these calculations are just formal, since $\psi$ doesn't have the necessary regularity to implement them. 
\newline
Furthermore, it is well known that the energy
\begin{equation}\label{eq:en_schr}
E(t):=\int_{\real^3}\frac{\hbar^2}{2}|\nabla\psi|^2+\frac{2}{p+1}|\psi|^{p+1}\de x
+\mez\int_{\real^3_x\times\real^3_y}|\psi(t,y)|^2\1{|x-y|}|\psi(t,x)|^2\de x\de y
\end{equation}
is conserved for the Schr\"odinger-Poisson system \eqref{eq:NLS}. Thus, it only remains to note that by Lemma \ref{lemma:phase} the energy in \eqref{eq:en_schr} is equal that one in \eqref{eq:en}.
\begin{dimo}
Let us consider the Schr\"odinger-Poisson system \eqref{eq:NLS}.
From standard theory about nonlinear Schr\"odinger equations it is well known that \eqref{eq:NLS} is globally well-posed for initial data in the space of energy: $\psi(0)=\psi_0\in H^1(\real^3)$. Now, let us take a sequence of mollifiers $\{\chi_\eps\}$ converging to the Dirac mass and define 
$\psi^\eps:=\chi_\eps\star\psi$, where $\psi\in\cfun(\real; H^1(\real^3))$ is the solution to \eqref{eq:NLS}. Then $\psi^\eps\in\cfun^\infty(\real^{1+3})$ and moreover
\begin{equation}\label{eq:NLS_eps}
i\hbar\d_t\psi^\eps+\frac{\hbar^2}{2}\Delta\psi^\eps=\chi_\eps\star\tonde{|\psi|^{p-1}\psi+V\psi},
\end{equation}
where $V$ is the classical Hartree potential. Therefore differentiating $|\psi^\eps|^2$ with respect to time we get
\begin{multline*}
\d_t|\psi^\eps|^2=2\re(\conj{\psi^\eps}\d_t\psi^\eps)
=2\re\tonde{\conj{\psi^\eps}(\frac{i\hbar}{2}\Delta\psi^\eps-\frac{i}{\hbar}\chi_\eps\star(|\psi|^{p-1}\psi+V\psi))}\\
=-\hbar\diver\tonde{\im(\conj{\psi^\eps}\nabla\psi^\eps)}
+\frac{2}{\hbar}\im\tonde{\conj{\psi^\eps}\chi_\eps\star(|\psi|^{p-1}\psi+V\psi)}.
\end{multline*}
If we differentiate with respect to time $\hbar\im(\conj{\psi^\eps}\nabla\psi^\eps)$ we get
\begin{multline*}
\d_t\tonde{\hbar\im(\conj{\psi^\eps}\nabla\psi^\eps)}
=\hbar\im\quadre{\tonde{-\frac{i\hbar}{2}\Delta\conj{\psi^\eps}
+\frac{i}{\hbar}\chi_\eps\star(|\psi|^{p-1}\conj{\psi})
+\frac{i}{\hbar}\chi_\eps\star(V\psi)}}\\
+\hbar\im\quadre{\conj{\psi^\eps}\tonde{\frac{i\hbar}{2}\delta\nabla\psi^\eps
-\frac{i}{\hbar}\chi_\eps\star\nabla(|\psi|^{p-1}\psi)-\frac{i}{\hbar}\chi_\eps\star\nabla(V\psi)}}\\
=\frac{\hbar^2}{2}\re\tonde{-\nabla\psi^\eps\Delta\conj{\psi^\eps}+\conj{\psi^\eps}\Delta\nabla\psi^\eps}
+\re\tonde{\chi_\eps\star(|\psi|^{p-1}\conj{\psi})\nabla\psi^\eps
-\conj{\psi^\eps}\chi_\eps\star\nabla(|\psi|^{p-1}\psi)}\\
+\re\tonde{\chi\star(V\conj{\psi})\nabla\psi^\eps-\chi_\eps\star\nabla(V\psi)\conj{\psi^\eps}}
=:A+B+C
\end{multline*}
Let us discuss these three terms separately. For the first term $A$ it is immediate that
\begin{equation*}
\frac{\hbar^2}{2}\re\tonde{-\nabla\psi^\eps\Delta\conj{\psi^\eps}+\conj{\psi^\eps}\Delta\nabla\psi^\eps}
=\hbar^2\diver\tonde{\re(\nabla\psi^\eps\otimes\nabla\conj{\psi^\eps})}
+\frac{\hbar^2}{4}\Delta\nabla|\psi^\eps|^2.
\end{equation*}
Let us discuss the second and the third terms, $B$ and $C$. We can recast $B$ in the following way
\begin{multline*}
\chi_\eps\star(|\psi|^{p-1}\conj{\psi})\nabla\psi^\eps
-\conj{\psi^\eps}\chi_\eps\star\nabla(|\psi|^{p-1}\psi)\\
=|\psi^\eps|^{p-1}\conj{\psi^\eps}\nabla\psi^\eps
-\conj{\psi^\eps}\chi_\eps\star\tonde{\psi\nabla|\psi|^{p-1}+|\psi|^{p-1}\nabla\psi}+R^\eps_1\\
=-\conj{\psi^\eps}\chi_\eps\star\tonde{\psi\nabla|\psi|^{p-1}}+R^\eps_1+R^\eps_2
=-|\psi^\eps|^{2}\nabla|\psi^\eps|^{p-1}+R^\eps_1+R^\eps_2+R^\eps_3,
\end{multline*}
where
\begin{gather*}
R^\eps_1:=\tonde{\chi_\eps\star(|\psi|^{p-1}\conj{\psi})-|\psi^\eps|^{p-1}\conj{\psi^\eps}}\nabla\psi^\eps\\
R^\eps_2:=\conj{\psi^\eps}\tonde{|\psi^\eps|^{p-1}\nabla\psi^\eps
-\chi_\eps\star(|\psi|^{p-1}\nabla\psi)}\\
R^\eps_3:=\conj{\psi^\eps}\tonde{\psi^\eps\nabla|\psi^\eps|^{p-1}-\chi_\eps\star(\psi\nabla|\psi|^{p-1})}.
\end{gather*}
Now it only remains to analyze how the remainder terms $R^\eps_j$ go to zero. We will prove it in the next Lemma. Furthermore, the following identity holds
\begin{equation*}
|\psi^\eps|^2\nabla|\psi^\eps|^{p-1}=\frac{p-1}{p+1}\nabla|\psi^\eps|^{p+1}.
\end{equation*}
For the third term $C$ after similar computations we get
\begin{equation*}
C=|\psi^\eps|\nabla(\chi_\eps\star V)+R^\eps_4,
\end{equation*}
where the remainder term $R^\eps_4$ will be analyzed in the next Lemma. Hence we can conclude that $\hbar\im(\conj{\psi^\eps}\nabla\psi^\eps)$ satisfies
\begin{multline*}
\d_t\tonde{\hbar\im(\conj{\psi^\eps}\nabla\psi^\eps)}
=-\hbar^2\diver\tonde{\re(\nabla\psi^\eps\otimes\nabla\conj{\psi^\eps})}
+\frac{\hbar^2}{4}\Delta\nabla|\psi^\eps|^2\\+\frac{p-1}{p+1}\nabla|\psi^\eps|^{p+1}
+|\psi^\eps|^2\nabla\tonde{\chi_\eps\star V}+R^\eps,
\end{multline*}
where $R^\eps:=R^\eps_1+R^\eps_2+R^\eps_3+R^\eps_4$. Now, as $\eps\to0$, we get
\begin{multline*}
\d_t\tonde{\hbar\im(\conj{\psi}\nabla\psi)}
=-\hbar^2\diver\tonde{\re(\nabla\psi\otimes\nabla\conj{\psi})}
+\frac{\hbar^2}{4}\Delta\nabla|\psi|^2\\+\frac{p-1}{p+1}\nabla|\psi|^{p+1}
+|\psi|^2\nabla V.
\end{multline*}
This identity is equivalent to \eqref{eq:curr_dens} since by the Lemma \ref{lemma:phase} we have
\begin{equation*}
\hbar^2\re\tonde{\nabla\psi\otimes\nabla\psi}=\hbar^2\nabla\sqrt{\rho}\otimes\nabla\sqrt{\rho}
+\Lambda\otimes\Lambda,
\end{equation*}
where $\Lambda$ is defined as in the Lemma \ref{lemma:phase} and 
$\sqrt{\rho}\Lambda=J$, where $J:=\hbar\im(\conj{\psi}\nabla\psi)$. 
\newline
Finally, let us note that, by the definition of $J$, we have
\begin{equation*}
\nabla\wedge J=\hbar\im(\nabla\conj{\psi}\wedge\nabla\psi)
\end{equation*}
then by the Corollary \ref{cor:irrot}, we get $\nabla\wedge J=2\nabla\sqrt{\rho}\wedge\Lambda$. 
\end{dimo}
\begin{lemma}
Let $0<T<\infty$, then
\begin{equation*}
\|R^\eps_j\|_{L^1_{t,x}([0, T]\times\real^3)}\to0,\qquad\textrm{as}\;\eps\to0,
\end{equation*}
$j=1, 2, 3, 4$.
\end{lemma}
\begin{dimo}
It is a direct consequence of Strichartz estimates for the solution of the Schr\"odinger-Poisson system \eqref{eq:NLS}
\begin{equation}
\|\psi\|_{L^q_tW^{1, r}([0, T]\times\real^3)}\leq C(E_0, \|\psi_0\|_{L^2(\real^3)}, T),
\end{equation}
where $(q, r)$ is a pair of admissible exponents and $E_0$ is the initial energy for \eqref{eq:NLS}. The first error can be controlled in the following way
\begin{multline*}
\|R^\eps_1\|_{L^1_{t, x}([0, T]\times\real^3)}\\
\leq
\lplq{\chi_\eps\star(|\psi|^{p-1}\conj{\psi})-|\psi^\eps|^{p-1}\conj{\psi^\eps}}{\frac{4(p+1)}{p+7}}{\frac{p+1}{p}}{[0, T]}
\lplq{\nabla\psi^\eps}{\frac{4(p+1)}{3(p-1)}}{p+1}{[0, T]}\\
\leq T^{\frac{p+7}{4(p+1)}}
\bigg(\lplqs{\chi_\eps\star(|\psi|^{p-1}\conj{\psi})-|\psi|^{p-1}\conj{\psi}}{\infty}{p+1}\\+
\lplqs{|\psi|^{p-1}\conj{\psi}-|\psi^\eps|^{p-1}\conj{\psi^\eps}}{\infty}{p+1}\bigg)
\lplqs{\nabla\psi^\eps}{\frac{3(p-1)}{4(p+1)}}{p+1}.
\end{multline*}
The remaining error terms can be computed in the same way. We remark that by using the Strichartz estimates it folllows 
$\psi\nabla|\psi|^{p-1}$ lies in $L^{\frac{4(p+1)}{3(p-1)}}_tL^{p+1}_x([0, T]\times\real^3)$.
\end{dimo}
\section{the fractional step: definitions and consistency}\label{sect:fract}
In this section we make use of the results of the previous sections to construct a sequence of approximate solutions of the QHD system
\begin{equation}\label{eq:QHD}
\left\{\begin{array}{l}
\d_t\rho+\diver J=0\\
\d_J+\diver\tonde{\frac{J\otimes J}{\rho}}+\nabla P(\rho)+\rho\nabla V+J=
\frac{\hbar^2}{2}\rho\nabla\tonde{\frac{\Delta\sqrt{\rho}}{\sqrt{\rho}}}\\
-\Delta V=\rho
\end{array}\right.
\end{equation}
with Cauchy data
\begin{equation}\label{eq:QHD_IV}
(\rho(0), J(0))=(\rho_0, J_0).
\end{equation}
\begin{defn}
Let $0<T<\infty$. 
We say  $\{(\rho^\tau, J^\tau)\}$ is a \emph{sequence of approximate solutions} to the system \eqref{eq:QHD3} in $[0, T)\times\real^3$, with initial data $(\rho_0, J_0)\in L^1_{loc}(\real^3)$, if there exist locally integrable functions 
$\sqrt{\rho^\tau}, \Lambda^\tau$, such that
$\sqrt{\rho^\tau}\in L^2_{loc}([0, T);H^1_{loc}(\real^3)), \Lambda^\tau\in L^2_{loc}([0, T);L^2_{loc}(\real^3))$ and if we define $\rho^\tau:=(\sqrt{\rho^\tau})^2$, $J^\tau:=\sqrt{\rho^\tau}\Lambda^\tau$, then
\begin{itemize}
\item for any test function $\eta\in\cfun_0^\infty([0,T)\times\real^3)$ one has
 \begin{equation}\label{eq:approx_QHD1}
 \int_0^T\int_{\real^3}\rho^\tau\d_t\eta+J^\tau\cdot\nabla\eta\de x\de t
 +\int_{\real^3}\rho_0\eta(0)\de x=o(1)
 \end{equation}
 as $\tau\to0$,
 \item for any test function $\zeta\in\cfun_0^\infty([0,T)\times\real^3;\real^3)$ we have
 \begin{multline}\label{eq:approx_QHD2}
 \int_0^T\int_{\real^3}J^\tau\cdot\d_t\zeta+\Lambda^\tau\otimes\Lambda^\tau:\nabla\zeta
 +P(\rho^\tau)\diver\zeta
 -\rho^\tau\nabla V^\tau\cdot\zeta-J^\tau\cdot\zeta\\
 +\hbar^2\nabla\sqrt{\rho^\tau}\otimes\nabla\sqrt{\rho^\tau}:\nabla\zeta
 -\frac{\hbar^2}{4}\rho^\tau\Delta\diver\zeta\de x\de t+\int_{\real^3}J_0\cdot\zeta(0)\de x=o(1)
 \end{multline}
 as $\tau\to0$;
 \item \emph{(generalized irrotationality condition)} for almost every $t\in(0, T)$ we have
 \begin{equation*}
 \nabla\wedge J^\tau=2\nabla\sqrt{\rho^\tau}\wedge\Lambda^\tau.
 \end{equation*}
 \end{itemize}
\end{defn}
Our fractional step method is based on the following simple idea. We split the evolution of our problem into two separate steps. Let us fix a (small) parameter $\tau>0$ then in the former step we solve a non-collisional QHD problem, while in the latter one we solve the collisional problem without QHD, and at this point we can start again with the non-collisional QHD problem, the main difficulty being the updating of the initial data at each time step. Indeed, as remarked in the previous section we are able to solve the non-collisional QHD only in the case of Cauchy data compatible with the Schr\"odinger picture. This restriction impose to reconstruct a wave function at each time step.
\newline
Now we explain how to set up the fractional step procedure which generates the approximate solutions. We first remark that, as in the previous section, this method requires special type of initial data $(\rho_0, J_0)$, namely we assume there exists $\psi_0\in H^1(\real^3)$, such that the hydrodynamic initial data are given via the Madelung transforms
\begin{equation*}
\rho_0=|\psi_0|^2,\qquad J_0=\hbar\im(\conj{\psi_0}\nabla\psi_0).
\end{equation*}
This assumption is physically relevant since it implies the compatibility of our solutions to the QHD problem with the wave mechanics approach. The iteration procedure can be defined in this following way. First of all, we take $\tau>0$, which will be the time mesh unit; therefore we define the approximate solutions in each time interval $[k\tau, (k+1)\tau)$, for any integer $k\geq0$.
\newline
At the first step, $k=0$, we solve the Cauchy problem for the Schr\"odinger-Poisson system
\begin{equation}\label{eq:NLS2}
\left\{\begin{array}{l}
i\hbar\d_t\psi^\tau+\frac{\hbar^2}{2}\Delta\psi^\tau=|\psi^\tau|^{p-1}\psi^\tau+V^\tau\psi^\tau\\
-\Delta V^\tau=|\psi^\tau|^2\\
\psi^\tau(0)=\psi_0
\end{array}\right.
\end{equation}
by looking for the restriction in $[0, \tau)$ of the unique strong solution 
$\psi\in\cfun(\real; H^1(\real^3))$ (see \cite{Caz}).
Let us define 
$\rho^\tau:=|\psi^\tau|^2, J^\tau:=\hbar\im(\conj{\psi^\tau}\nabla\psi^\tau)$. Then, as shown in the previous section, $(\rho^\tau, J^\tau)$ is a weak solution to the non-collisional QHD system. 
Let us assume that we know $\psi^\tau$ in the space-time slab 
$[(k-1)\tau, k\tau)\times\real^3$, we want to set up a recursive method, hence we have to show how to define $\psi^\tau, \rho^\tau, J^\tau$ in the strip $[k\tau, (k+1)\tau)$.
\newline
In order to take into account the presence of the collisional term $f=J$ we update $\psi$ in $t=k\tau$, namely we define $\psi^\tau(k\tau+)$.
The construction of $\psi^\tau(k\tau+)$ will be done by means of the polar decomposition described in the Section \ref{sect:polar}.
\newline
Let us apply the Lemma \ref{lemma:updat}, with $\psi=\psi^\tau(k\tau-)$, 
$\eps=\tau2^{-k}\|\psi_0\|_{H^1(\real^3)}$, then we can define
\begin{equation}\label{eq:125}
\psi^\tau(k\tau+)=\tilde\psi,
\end{equation}
by using the wave function $\tilde\psi$ defined in the Lemma \ref{lemma:updat}. Therefore we have
\begin{align}
\rho^\tau(k\tau+)&=\rho^\tau(k\tau-)\\
\Lambda^\tau(k\tau+)&=(1-\tau)\Lambda^\tau(k\tau-)+R_k\label{eq:126},
\end{align}
where $\|R_k\|_{L^2(\real^3)}\leq\tau2^{-k}\|\psi_0\|_{H^1(\real^3)}$ and
\begin{equation}
\nabla\psi^\tau(k\tau+)=\nabla\psi^\tau(k\tau-)-i\frac{\tau}{\hbar}\phi^\ast\Lambda^\tau(k\tau-)+r_{k, \tau},
\end{equation}
with $\|\phi^\ast\|_{L^\infty}\leq1$ and
\begin{equation*}
 \|r_{k, \tau}\|_{L^2}\leq C(\tau\|\nabla\psi^\tau(k\tau-)\|+\tau2^{-k}\|\psi_0\|_{H^1(\real^3)})
 \lesssim\tau E_0^{\mez}.
 \end{equation*}
\newline
We then solve the Schr\"odinger-Poisson system with initial data $\psi(0)=\psi^\tau(k\tau+)$. We define $\psi^\tau$ in the time strip $[k\tau, (k+1)\tau)$ as the restriction of the Schr\"odinger-Poisson solution just found in $[0, \tau)$, furthermore, we define $\rho^\tau:=|\psi^\tau|^2, 
J^\tau:=\hbar\im(\conj{\psi^\tau}\nabla\psi^\tau)$ as the solution of the non-collisional QHD 
\eqref{eq:QHD_farl}.
\newline
With this procedure  we can go on every time strip and then construct an approximate solution 
$(\rho^\tau, J^\tau, V^\tau)$ of the QHD system.
\begin{oss}\label{oss:gen_case}
To cover the general case where
$f(\sqrt{\rho}, J, \nabla\sqrt{\rho})=\alpha J+\rho\nabla g(t, x, \sqrt{\rho}, \Lambda, \nabla\sqrt{\rho})$,
$\alpha\geq0$, $g$ satisfies the conditions in the Remark \ref{oss:coll} and one has a non-zero doping profile, we solve the following Cauchy problem
\begin{equation}\label{eq:NLS3}
\left\{\begin{array}{l}
i\hbar\d_t\psi+\frac{\hbar^2}{2}\Delta\psi=|\psi|^{p-1}\psi+V\psi+g\psi\\
-\Delta V=|\psi|^2-C(x)\\
\psi(0)=\psi_0
\end{array}\right.
\end{equation}
The global well-posedness follows from the result in \cite{Caz}.
\end{oss}
\begin{thm}[Consistency of the approximate solutions]\label{thm:cons}
Let us consider the sequence of approximate solutions $\{(\rho^\tau, J^\tau)\}_{\tau>0}$ constructed via the fractional step method, and assume there exists $0<T<\infty$, 
$\sqrt{\rho}\in L^2_{loc}([0, T);H^1_{loc}(\real^3))$ and 
$\Lambda\in L^2_{loc}([0, T);L^2_{loc}(\real^3))$, such that
\begin{align}
\sqrt{\rho^\tau}&\to\sqrt{\rho}\qquad\textrm{in}\;L^2([0, T);H^1_{loc}(\real^3))\\
\Lambda^\tau&\to\Lambda\qquad\textrm{in}\;L^2([0, T); L^2_{loc}(\real^3)).
\end{align}
Then the limit function $(\rho, J)$, where as before $J=\sqrt{\rho}\Lambda$, is a weak solution of the QHD system, with Cauchy data 
$(\rho_0, J_0)$.
\end{thm}
\begin{dimo}
Let us plug the approximate solutions $(\rho^\tau, J^\tau)$ in the weak formulation and let 
$\zeta\in\cfun_0^\infty([0, T)\times\real^3)$, then we have
\begin{gather*}
\begin{split}
\int_0^\infty\int_{\real^3}J^\tau\cdot\d_t\zeta+\Lambda^\tau\otimes\Lambda^\tau:\nabla\zeta
+P(\rho^\tau)\diver\zeta-\rho^\tau\nabla V^\tau\cdot\zeta-J^\tau\cdot\zeta\\
+\hbar^2\nabla\sqrt{\rho^\tau}\otimes\nabla\sqrt{\rho^\tau}:\nabla\zeta
-\frac{\hbar^2}{4}\rho^\tau\Delta\diver\zeta\de x\de t+\int_{\real^3}J_0\cdot\zeta(0)\de x
\end{split}\\
\begin{split}
=\sum_{k=0}^\infty\slabint J^\tau\cdot\d_t\zeta+\Lambda^\tau\otimes\Lambda^\tau:\nabla\zeta
+P(\rho^\tau)\diver\zeta-\rho^\tau\nabla V^\tau\cdot\zeta-J^\tau\cdot\zeta\\
+\hbar^2\nabla\sqrt{\rho^\tau}\otimes\nabla\sqrt{\rho^\tau}:\nabla\zeta
-\frac{\hbar^2}{4}\rho^\tau\Delta\diver\zeta\de x\de t
+\int_{\real^3}J_0\cdot\zeta(0)\de x
\end{split}\\
\begin{split}
=\sum_{k=0}^\infty\bigg[\slabint-J^\tau\cdot\zeta\de x\de t
+\int_{\real^3}J^\tau((k+1)\tau-)\cdot\zeta((k+1)\tau)\\
-J^\tau(k\tau)\cdot\zeta(k\tau)\de x\bigg]
+\int_{\real^3}J_0\cdot\zeta(0)\de x
\end{split}\\
\end{gather*}
\begin{gather*}
\begin{split}
=\sum_{k=0}^\infty\slabint-J^\tau\cdot\zeta\de x\de t
+\sum_{k=1}^\infty\int_{\real^3}\tonde{J^\tau(k\tau-)-J^\tau(k\tau)}\cdot\zeta(k\tau)\de x
\end{split}\\
\begin{split}
=\sum_{k=0}^\infty\slabint-J^\tau\cdot\zeta\de x\de t
+\sum_{k=1}^\infty\tau\int_{\real^3}J^\tau(k\tau-)\cdot\zeta(k\tau)\de x\\
+\sum_{k=1}^\infty\int_{\real^3}(1-\tau)\sqrt{\rho^\tau(k\tau)}R_k\cdot\zeta(k\tau)\de x
\end{split}\\
\begin{split}
=\sum_{k=0}^\infty\slabint J^\tau((k+1)\tau-)\cdot\zeta((k+1)\tau)-J^\tau(t)\cdot\zeta(t)\de x\de t\\
+(1-\tau)\sum_{k=0}^\infty\int_{\real^3}\sqrt{\rho^\tau}(k\tau)R_k\cdot\zeta(k\tau)\de x
\end{split}\\
=o(1)+O(\tau)
\end{gather*}
as $\tau\to0$.
\end{dimo}
\section{a priori estimates and convergence}\label{sect:apr}
In this section we obtain various a priori estimates necessary to show the compactness of the sequence of approximate solutions $(\rho^\tau, J^\tau)$ in some appropriate function spaces. As we stated in the Theorem \ref{thm:cons}, we wish to prove the strong convergence of $\{\sqrt{\rho^\tau}\}$ in $L^2_{loc}([0, T);H^1_{loc}(\real^3))$ and that one of 
$\{\Lambda^\tau\}$ in $L^2_{loc}([0, T);L^2_{loc}(\real^3))$. To achieve this goal we use a compactness result in the class of the Aubin-Lions's type lemma, due to Rakotoson-Temam \cite{RT} (see the Section 2). The plan of this section is the following, first of all we get a discrete version of the (dissipative) energy inequality for the system \eqref{eq:QHD3}, later we use the Strichartz estimates for $\nabla\psi^\tau$ by means of the formula \eqref{eq:131} below. Consequently via the Strichartz estimates and by using the local smoothing results of Theorems \ref{thm:smooth1}, \ref{thm:smooth2}, we deduce some further regularity properties of the sequence $\{\nabla\psi^\tau\}$. This fact will be stated in the next Proposition 
\ref{prop:28}. In this way it is possible to get the regularity properties of the sequence $\{\nabla\psi^\tau\}$ which are needed to apply the Theorem \ref{thm:comp} and hence to get the convergence for $\{\sqrt{\rho^\tau}$ and $\{\Lambda^\tau\}$.
\newline
Let us begin with the energy inequality. First of all note that if we have a sufficiently regular solution of the QHD system, one has
\begin{equation}
E(t)=-\int_0^t\int_{\real^3}|\Lambda|^2\de x\de t'+E_0,
\end{equation}
where the energy is defined as in \eqref{eq:en}.
Now we would like to find a discrete version of the energy dissipation for the approximate solutions. 
\begin{prop}[Discrete energy inequality]
Let $(\rho^\tau, J^\tau)$ be an approximate solution of the QHD system, with $0<\tau<1$. Then, for 
$t\in[N\tau, (N+1)\tau)$ we have
\begin{equation}\label{eq:en-diss}
E^\tau(t)\leq-\frac{\tau}{2}\sum_{k=1}^N\|\Lambda(k\tau-)\|_{L^2(\real^3)}+(1+\tau)E_0.
\end{equation}
\end{prop}
\begin{dimo}
For all $k\geq1$, we have
\begin{align*}
E^\tau(k\tau+)-E^\tau(k\tau-)=&
\int\mez|\Lambda^\tau(k\tau+)|^2-\mez|\Lambda^\tau(k\tau-)|^2\de x\\
=&\mez\int(-2\tau+\tau^2)|\Lambda^\tau(k\tau-)|^2\\
&\qquad+2(1-\tau)\Lambda^\tau(k\tau-)\cdot R_k+|R_k|^2\de x\\
\leq&\mez\int(-2\tau+\tau^2)|\Lambda^\tau(k\tau-)|^2\de x+(1-\tau)\alpha|\Lambda^\tau(k\tau-)|^2\\
&\qquad+\frac{1-\tau}{\alpha}|R_k|^2+|R_k|^2\de x\\
=&\mez\int(-2\tau+\tau^2+\alpha-\alpha\tau)|\Lambda^\tau(k\tau-)|^2\\
&\qquad+\tonde{\frac{1-\tau+\alpha}{\alpha}}|R_k|^2\de x.
\end{align*}
Here $R_k$ denotes the error term as in espression \eqref{eq:126}.
If we choose $\alpha=\tau$, it follows
\begin{multline}\label{eq:124}
E^\tau(k\tau+)-E^\tau(k\tau-)\leq-\frac{\tau}{2}\|\Lambda^\tau(k\tau-)\|_{L^2}^2
+\1{2\tau}\|R_k\|^2_{L^2}\\
\leq-\frac{\tau}{2}\|\Lambda^\tau(k\tau-)\|^2_{L^2}+\tau2^{-k-1}\|\psi_0\|_{H^1}.
\end{multline}
The inequality \eqref{eq:en-diss} follows by summing up all the terms in \eqref{eq:124} and by the energy conservation in each time strip $[k\tau, (k+1)\tau)$.
\end{dimo}
Unfortunately the energy estimates are not sufficient to get enough compactness to show the convergence of the sequence of the approximate solutions. Indeed from the discrete energy inequality, we get only the weak convergence of $\nabla\psi^\tau$ in $L^\infty([0, \infty);H^1(\real^3))$, and therefore the quadratic terms in \eqref{eq:approx_QHD2} could exhibit some concentrations phenomena in the limit.
\newline
More precisely, 
from energy inequality we get the sequence $\{\psi^\tau\}$ is uniformly bounded in 
$L^\infty([0,\infty);H^1(\real^3))$, hence there exists $\psi\in L^\infty([0, \infty);H^1(\real^3))$, such that
\begin{equation*}
\psi^\tau\wlim\psi\qquad\textrm{in}\;L^\infty([0,\infty);H^1(\real^3)).
\end{equation*}
Therefore we get
\begin{align*}
\sqrt{\rho^\tau}&\wlim\sqrt{\rho}\qquad\textrm{in}\;L^\infty([0,\infty);H^1(\real^3))\\
\Lambda^\tau&\wlim\Lambda\qquad\textrm{in}\;L^\infty([0,\infty);L^2(\real^3))\\
J^\tau&\wlim\sqrt{\rho}\Lambda\qquad\textrm{in}\;L^\infty([0, \infty);L^1_{loc}(\real^3)).
\end{align*}
The need to pass into the limit the quadratic expressions leads us to look for a priori estimates in stronger norms. The relationships with the Schr\"odinger equation brings naturally into this search the Strichartz-type estimates. The following results are concerned with these estimates. However they are \emph{not} an immediate consequence of the Strichartz estimates for the Schr\"odinger equation since we have to take into account the effects of the updating procedure which we implement at each time step.
\newline
The following Remark \ref{rmk:25} will clarify in detail why we stated the Lemma \ref{lemma:updat} in place of using the exact factorization property provided by the Lemma \ref{lemma:phase}.
\begin{lemma}\label{lemma:28}
Let $\psi^\tau$ be the wave function defined by the fractional step method (see Section 
\ref{sect:fract}), and let $t\in[N\tau, (N+1)\tau)$. Then we have
\begin{multline}\label{eq:131}
\nabla\psi^\tau(t)=U(t)\nabla\psi_0
-i\frac{\tau}{\hbar}\sum_{k=1}^NU(t-k\tau)\tonde{\phi^\tau_k\Lambda^\tau(k\tau-)}\\
-i\int_0^tU(t-s)F(s)\de s+\sum_{k=1}^NU(t-k\tau)r^\tau_k,
\end{multline}
where as before $U(t)$ is the free Schr\"odinger group, 
\begin{equation}
\|\phi^\tau_k\|_{L^\infty(\real^3)}\leq1, \quad\|r^\tau_k\|_{L^2(\real^3)}\leq\tau\|\psi_0\|_{H^1(\real^3)}
\end{equation}
and $F=\nabla(|\psi^\tau|^{p-1}\psi^\tau+V^\tau\psi^\tau)$.
\end{lemma}
\begin{oss}\label{oss:25}
Clearly, with a general collision term $f(\sqrt{\rho}, J, \nabla\sqrt{\rho})$ (see the Remark \ref{oss:coll}) and a non-zero doping profile $C(x)$ one has $F$ defined by 
$F=\nabla(|\psi^\tau|^{p-1}\psi^\tau+V^\tau\psi^\tau+g^\tau\psi^\tau)$, where now $V^\tau$ solves
\begin{equation*}
-\Delta V^\tau=|\psi^\tau|^2-C(x)
\end{equation*}
and $g^\tau=g(t, x, \sqrt{\rho^\tau}, \Lambda^\tau, \nabla\sqrt{\rho^\tau})$.
\end{oss}
\begin{oss}\label{rmk:25}
Now we can see why the updating step in \eqref{eq:125} has been defined through the Lemma 
\ref{lemma:updat}. 
Indeed, one could possibly try to use the exact unitary factor as in the Lemma \ref{lemma:phase}, namely 
$\psi^\tau(k\tau+):=\phi^{(1-\tau)}\sqrt{\rho}$, where $\phi, \sqrt{\rho}$ are respectively the unitary factor and the amplitude of $\psi^\tau(k\tau-)$ (note that $\phi$ is uniquely determined $\sqrt{\rho}\de x$, and 
$|\phi|=1$, $\sqrt{\rho}\de x$-a.e.) The meaning of this definition is clear if we approximate 
$\psi^\tau(k\tau-)$ (as in the proof of the Lemma \ref{lemma:updat}) with smooth 
$\psi_n=e^{i\theta_n}\sqrt{\rho_n}$ and then we update $\psi_n$ with 
$\tilde\psi_n:=e^{i(1-\tau)\theta_n}\sqrt{\rho_n}$. Hence $\psi(k\tau+):=\lim_{n\to\infty}\tilde\psi_n$ in 
$H^1(\real^3)$. Moreover in this case we get the exact updating formulas
\begin{align*}
|\psi^\tau(k\tau+)|^2&=|\psi^\tau(k\tau)|^2\\
\hbar\im(\conj{\psi^\tau}\nabla\psi^\tau)(k\tau+)&=(1-\tau)\hbar\im(\conj{\psi^\tau}\nabla\psi^\tau)(k\tau-),
\end{align*}
without the small error $r_\eps$ in the expression for the density current, however instead of the formula 
\eqref{eq:131} we would find
\begin{align*}
\nabla\psi^\tau(t)=&U(t-N\tau)\sigma^\tau_{N\tau}U(\tau)\dotsc\sigma^\tau_\tau U(\tau)\nabla\psi_0\\
&-i\frac{\tau}{\hbar}U(t-N\tau)\sigma^\tau_{N\tau}\Lambda^\tau(N\tau-)
+\dotsc-i\frac{\tau}{\hbar}U(t-N\tau)\dotsc U(\tau)\phi^\tau_\tau\Lambda^\tau(\tau-)\\
&-i\int_{N\tau}^tU(t-s)F(s)\de s
-iU(t-N\tau)\sigma^\tau_{N\tau}\int_{(N-1)\tau}^{N\tau}U(N\tau-s)F(s)\de s\\
&+\dotsc
-iU(t-N\tau)\sigma^\tau_{N\tau}U(\tau)\dotsc\sigma^\tau_\tau\int_0^\tau U(\tau-s)F(s)\de s,
\end{align*}
where $\sigma^\tau_{k\tau}=(\phi^\tau_{k\tau})^{-t}$, $\phi^\tau_{k\tau}$ being the polar factor of 
$\psi^\tau(k\tau-)$.
Now, to recover an expression similar to \eqref{eq:131}, we have to calculate the commutators between the free Schr\"odinger evolution operator and the multiplication operators by $\sigma^\tau_{k\tau}$. The estimates of these commutators are not known for non-smooth $\sigma^\tau_{k\tau}$.
\end{oss}
\begin{dimo}
Since $\psi^\tau$ is solution of the Schr\"odinger-Poisson system in the space-time slab 
$[N\tau, (N+1)\tau)\times\real^3$, then we can write
\begin{equation}\label{eq:111}
\nabla\psi^\tau(t)=U(t-N\tau)\nabla\psi^\tau(N\tau+)-i\int_{N\tau}^tU(t-s)F(s)\de s,
\end{equation}
where $F$ is defined in the statement of the Lemma \ref{lemma:28}. Now there exists a piecewise smooth function 
$\theta_N$, as specified in the proof of the Lemma \ref{lemma:updat}, such that
\begin{equation*}
\psi(N\tau+)=e^{i(1-\tau)\theta_N}\sqrt{\rho_n},
\end{equation*}
and furthermore all the estimates of the Lemma \ref{lemma:updat} hold, with $\psi=\psi^\tau(N\tau-)$, 
$\tilde\psi=\psi^\tau(N\tau+)$ and $\eps=2^{-N}\tau\|\psi_0\|_{H^1(\real^3)}$. Therefore, we have
\begin{equation}\label{eq:112}
\nabla\psi^\tau(N\tau+)=\nabla\psi^\tau(N\tau-)
-i\frac{\tau}{\hbar}e^{i(1-\tau)\theta_N}\Lambda^\tau(N\tau-)+r^\tau_N,
\end{equation}
where $\|r^\tau_N\|_{L^2}\leq\tau\|\psi_0\|_{H^1}$.
By plugging \eqref{eq:112} into \eqref{eq:111} we deduce
\begin{multline*}
\nabla\psi^\tau(t)=U(t-N\tau)\nabla\psi^\tau(N\tau-)
-i\frac{\tau}{\hbar}U(t-N\tau)(e^{i(1-\tau)\theta_N}\Lambda^\tau(N\tau-))\\
+U(t-N\tau)r^\tau_N-i\int_{N\tau}^tU(t-s)F(s)\de s.
\end{multline*}
Let us iterate this formula, repeating the same procedure for $\nabla\psi^\tau(N\tau-)$, then \eqref{eq:131} holds.
\end{dimo}
At this point we can use the formula \eqref{eq:131} to obtain Strichartz estimates for $\nabla\psi^\tau$, simply by applying the Theorem \ref{thm:strich} to each of the terms in \eqref{eq:131}. After some computations, the following result holds.
\begin{prop}[Strichartz estimates for $\nabla\psi^\tau$]\label{prop:26}
Let, $0<T<\infty$, let $\psi^\tau$ be as in the previous section, then one has
\begin{equation}
\lplq{\nabla\psi^\tau}{q}{r}{[0,T]}\leq C(E_0^{\mez}, \|\rho_0\|_{L^1(\real^3)}, T)
\end{equation}
for each admissible pair of exponents $(q,r)$.
\end{prop}
\begin{oss}\label{oss:strich_gen}
As we showed in Remark \ref{oss:25}, in the general case the non-homogeneous term $F$ is slightly different. Anyway by the Strichartz estimates the same result of Proposition \ref{prop:26} holds even in the case of modified non-homogeneous term $F$.
\end{oss}
\begin{dimo}
First of all, let us prove the Proposition \ref{prop:26}, for a small time $0<T_1\leq T$ and let 
$(q, r)$ be an admissible pair of exponents. We choose $T_1>0$ later. Let $N$ be a positive integer such that $T_1\leq N\tau$. By applying the Theorem \ref{thm:strich} to the formula \eqref{eq:131}, we get
\begin{align*}
\LpLq{\nabla\psi^\tau}{q}{r}{[0,T_1]}\leq&\LpLq{U(t)\nabla\psi_0}{q}{r}{[0,T_1]}\\
&+\frac{\tau}{\hbar}\sum_{k=1}^N
\LpLq{U(t-k\tau)\tonde{e^{i(1-\tau)\theta_k}\Lambda^\tau(k\tau-)}}{q}{r}{[0,T_1]}\\
&+\sum_{k=1}^N\LpLq{U(t-k\tau)r^\tau_k}{q}{r}{[0,T_1]}\\
&+\LpLq{\int_0^tU(t-s)F(s)\de s}{q}{r}{[0,T_1]}\\
=:A+B+C+D.
\end{align*}
Now we estimate term by term the above expression.
\newline
The estimate of $A$ is straightforward, since
\begin{equation*}
\|U(t)\nabla\psi_0\|_{L^q_tL^r_x([0,T_1]\times\real^3)}\lesssim\|\nabla\psi_0\|_{L^2(\real^3)}.
\end{equation*}
The estimate of $B$ follows from
\begin{multline}
\frac{\tau}{\hbar}\sum_{k=1}^N
\LpLq{U(t-k\tau)\tonde{e^{i(1-\tau)\theta_k}\Lambda^\tau(k\tau-)}}{q}{r}{[0,T_1]}\\
\lesssim\tau\sum_{k=1}^N\|\Lambda^\tau(k\tau-)\|_{L^2(\real^3)}\lesssim T_1E_0^{\mez}.
\end{multline}
The term $C$ can be estimated in a similar way, namely
\begin{equation}
\sum_{k=1}^N\|U(t-k\tau)r^\tau_k\|_{L^q_tL^r_x}\lesssim\sum_{k=1}^N\|r^\tau_k\|_{L^2(\real^3)}
\lesssim T_1\|\psi_0\|_{H^1(\real^3)}.
\end{equation}
The last term is a little bit more tricky to estimate. First of all we decompose $F$ into three terms,
$F=F_1+F_2+F_3$, where $F_1=\nabla(|\psi^\tau|^{p-1}\psi^\tau)$, $F_2=\nabla V^\tau\psi^\tau$ and
$F_3=V^\tau\nabla\psi^\tau$.
By the Strichartz estimates (Theorem \ref{thm:strich}), we have
\begin{multline}
\LpLq{\int_0^tU(t-s)F(s)\de s}{q}{r}{[0,T_1]}\\
\lesssim\LpLq{F_1}{q_1'}{r_1'}{[0,T_1]}+\LpLq{F_2}{q_2'}{r_2'}{[0,T_1]}+\LpLq{F_3}{q_3'}{r_3'}{[0,T_1]},
\end{multline}
where $(q_i, r_i)$ are pairs of admissible exponents.
Let us start with the first term.
\begin{lemma}
There exists $\alpha>0$, depending on $p$, such that
\begin{equation}\label{eq:141}
\lplq{|\psi^\tau|^{p-1}\nabla\psi^\tau}{\tilde q'}{\tilde r'}{[0,T_1]}\lesssim T_1^\alpha
\|\psi^\tau\|_{\dot S^1([0, T_1]\times\real^3)}.
\end{equation}
\end{lemma}
\begin{dimo}
First of all let us apply H\"older inequality on the left-hand side of \eqref{eq:141}. Then we have
\begin{multline}\label{eq:201}
\lplq{|\psi^\tau|^{p-1}\nabla\psi^\tau}{\tilde q'}{\tilde r'}{[0,T_1]}\\\lesssim T_1^\alpha
\lplq{|\psi^\tau|^{p-1}}{q_1}{r_1}{[0,T_1]}\lplq{\nabla\psi}{q_2}{r_2}{[0,T_1]}\\
=T_1^\alpha\lplq{\psi^\tau}{q_1(p-1)}{r_1(p-1)}{[0,T_1]}^{p-1}
\lplq{\nabla\psi^\tau}{q_2}{r_2}{[0,T_1]}.
\end{multline}
Now we want $\1{q_1(p-1)}=\frac{3}{2}\tonde{\1{6}-\1{r_1(p-1)}}$ and 
$\1{q_2}=\frac{3}{2}\tonde{\mez-\1{r_2}}$, in such a way that 
$\lplqs{f}{q_1(p-1)}{r_1(p-1)}, \lplqs{\nabla f}{q_2}{r_2} \leq\|f\|_{\dot S^1}=\|\nabla f\|_{\dot S^0}$.
We already know 
$\1{\tilde q'}=1+\frac{3}{2}\tonde{\mez-\1{\tilde r'}}$, then putting together the conditions on $(\tilde q, \tilde r), (q_j, r_j)$, it follows
\begin{multline*}
\1{\tilde q'}=\1{\alpha}+\1{q_1}+\1{q_2}=\1{\alpha}+(p-1)\frac{3}{2}\tonde{\1{6}-\1{r_1(p-1)}}\\
+\frac{3}{2}\tonde{\mez-\1{r_2}}=1+\frac{3}{2}\tonde{\mez-\1{\tilde r'}}
\end{multline*}
and when $1\leq p<5$,
\begin{equation*}
\alpha=\frac{5-p}{4}>0.
\end{equation*}
This means that we can always choose pairs $(\tilde q', \tilde r'), (q_1, r_1), (q_2, r_2)$ satisfying the previous conditions so that the inequality \eqref{eq:201} holds, with $\alpha>0$.
\newline
For instance, if $1\leq p\leq3$, we can choose $\1{r_1}=\frac{p-1}{6}, \1{r_2}=\mez$, therefore 
$q_1=q_2=\infty$, hence we have $\1{\tilde r'}=\frac{2+p}{6}, \1{\tilde q'}=\frac{5-p}{4}$. In the case 
$3\leq p<5$, we take $\1{r_1}=\frac{p-1}{6}, \1{r_2}=\1{6}$ then $q_1=\infty, q_2=2$, hence we have $\1{\tilde r'}=\frac{p}{6}, \1{\tilde q'}=\frac{7-p}{4}$.
\end{dimo}
Now let us consider the second term: 
$V^\tau\nabla\psi^\tau$, here we choose $(q_2',r_2')=(1,2)$, so by the H\"older and the Hardy-Littlewood-Sobolev (see \cite{Sog}) inequalities one has
\begin{multline*}
\lplq{V^\tau\nabla\psi^\tau}{1}{2}{[0,T_1]}\leq T_1^{\mez}\lplqs{V^\tau}{\infty}{2}
\lplqs{\nabla\psi^\tau}{2}{6}\\
\lesssim T_1^{\mez}\lplqs{\psi^\tau}{\infty}{2}^2\lplqs{\nabla\psi^\tau}{2}{6}.
\end{multline*}
For the third term, we choose $(q_3', r_3')=(\frac{2}{2-3\eps}, \frac{2}{1+2\eps})$ and again by using the Hardy-Littlewood-Sobolev and the H\"older inequalities, we have
\begin{multline*}
\lplq{\nabla V^\tau\psi^\tau}{\frac{2}{2-3\eps}}{\frac{2}{1+2\eps}}{[0,T_1]}
\leq T_1^{\mez}\lplqs{\nabla V^\tau}{\frac{2}{1-3\eps}}{\1{\eps}}\lplqs{\psi^\tau}{\infty}{2}\\
\lesssim T_1^{\mez}\lplqs{\nabla|\psi^\tau|^2}{\frac{2}{1-3\eps}}{\frac{3}{2+3\eps}}
\lplqs{\tilde\psi^\tau}{\infty}{2}
\lesssim T_1^{\mez}\lplqs{\psi^\tau}{\infty}{2}^2
\lplqs{\nabla\psi^\tau}{\frac{2}{1-3\eps}}{\frac{6}{1+6\eps}}
\end{multline*}
Now, we summarize the previous estimates by using \eqref{eq:131} in the following way
\begin{multline}
\|\nabla\psi^\tau\|_{\dot S^0([0,T_1]\times\real^3)}\lesssim\|\nabla\psi_0\|_{L^2(\real^3)}+T_1E_0^{\mez}
+T_1^\alpha\|\nabla\psi^\tau\|_{\dot S^0([0,T_1]\times\real^3)}^p\\
+T_1^{\mez}\|\psi_0\|_{L^2(\real^3)}^2\|\nabla\psi^\tau\|_{\dot S^0([0,T_1]\times\real^3)}\\
\lesssim(1+T)E_0^{\mez}
+T_1^\alpha\|\nabla\psi^\tau\|_{\dot S^0([0,T_1]\times\real^3)}^p
+T_1^{\mez}\|\psi_0\|_{L^2(\real^3)}^2\|\nabla\psi^\tau\|_{\dot S^0([0,T_1]\times\real^3)}.
\end{multline}
\begin{lemma}\label{lemma:38}
There exist $T_1(E_0, \|\psi_0\|_{L^2(\real^3)}, T)>0$ and $C_1(E_0, \|\psi_0\|_{L^2(\real^3)}, T)>0$, independent on $\tau$, such that
\begin{equation}\label{eq:212}
\|\nabla\psi^\tau\|_{\dot S^0([0, \tilde T]\times\real^3)}\leq C_1(E_0, \|\psi_0\|_{L^2(\real^3)}, T)
\end{equation}
for all $0<\tilde T\leq T_1(E_0, \|\psi_0\|_{L^2(\real^3)})$.
\end{lemma}
Let us recall that $E^\tau(t)=E_0$ and $\|\psi^\tau(t)\|_{L^2(\real^3)}=\|\psi_0\|_{L^2(\real^3)}$, hence we are in the situation in which we can repeat our argument every time interval of length $T_1$, depending always on the same parameters $E_0, \|\psi_0\|_{L^2}$. The consequence of this fact is the following inequality on $[0, T]$
\begin{multline}\label{eq:213}
\|\nabla\psi^\tau\|_{\dot S^0([0, T]\times\real^3)}\\\leq C_1(E_0, \|\psi_0\|_{L^2}, T)
\tonde{\quadre{\frac{T}{T_1}}+1}=C(\|\psi_0\|_{L^2}, E_0, T).
\end{multline}
\end{dimo}
\emph{Proof of the Lemma \ref{lemma:38}}. Let us consider the non-trivial case $\|\psi_0\|_{L^2}>0$. Assume that $X\in(0, \infty)$ satisfies
\begin{equation}
X\leq A+\mu X+\lambda X^p=\phi(X),
\end{equation}
with $p>1$, $A>0$ and for all $0<\mu<1$, $\lambda>0$. Let $X_*$ be such that $\phi'(X_*)=1$, namely 
$X_*=\tonde{\frac{1-\mu}{p\lambda}}^{\1{p-1}}$, hence one has $\phi(X_*)<X_*$ each time the following inequality is satisfied
\begin{equation}\label{eq:231}
\tonde{\1{p^{\1{p-1}}}-\1{p^{\frac{p}{p-1}}}}\frac{(1-\mu)^{\frac{p}{p-1}}}{\lambda^{\1{p-1}}}>A.
\end{equation}
Therefore the convexity of $\phi$ implies that, if the condition \eqref{eq:231} holds, there exist two roots $X_{\pm}$, $X_+(\mu, \lambda, A)>X_*>X_-(\mu, \lambda, A)$, to the equation $\phi(X)=X$. It then follows either $0\leq X\leq X_-$, or $X\geq X_+$. In our case $\mu=T_1^{1/2}\|\psi_0\|_{L^2}^2$, 
$\lambda=T_1^\alpha$, $A=(1+T)E_0^{1/2}$, hence we assume
\begin{gather*}
\mu=T_1^{1/2}\|\psi_0\|_{L^2}^2<\mez\\
\lambda=T_1^\alpha=T_1^{\frac{5-p}{4}}
<\quadre{\frac{p^{-\1{p-1}}-p^{-\frac{p}{p-1}}}{2^{\frac{p}{p-1}}(1+T)E_0^{1/2}}}^{p-1}.
\end{gather*}
Therefore we choose
\begin{equation}
T_1:=\min\quadre{(2\|\psi_0\|_{L^2})^{-2}, 
\quadre{\frac{p^{-\1{p-1}}-p^{-\frac{p}{p-1}}}{2^{\frac{p}{p-1}}(1+T)E_0^{1/2}}}^{\frac{4(p-1)}{5-p}}}.
\end{equation}
Clearly we cannot have
\begin{equation*}
x_*=\quadre{\frac{1-T_1\|\psi_0\|_{L^2}}{pT_1^\alpha}}^{\1{p-1}}\leq x_+\leq
\|\nabla\psi^\tau\|_{\dot S^0([0, T_1]\times\real^3)},
\end{equation*}
since we get a contradiction as $T_1\to0$, hence
\begin{equation}
\|\nabla\psi^\tau\|_{\dot S^0([0, T_1]\times\real^3)}\leq X_-.
\end{equation}
\hfill $\square$\vspace{1cm}\newline
\begin{cor}
Let $0<T<\infty$ and let $\sqrt{\rho^\tau}, \Lambda^\tau$ be as in previous section, then 
\begin{equation}
\lplq{\nabla\sqrt{\rho^\tau}}{q}{r}{[0, T]}+\lplq{\Lambda^\tau}{q}{r}{[0, T]}
\leq C(E_0^{\mez}, \|\rho_0\|_{L^1(\real^3)}, T),
\end{equation}
for each admissible pair of exponents $(q, r)$.
\end{cor}
Unfortunately this is not enough to achieve the convergence of the quadratic terms. We need some additional compactness estimates on the sequence $\{\nabla\psi^\tau\}$ in order to apply Theorem 
\ref{thm:comp}. In particular we need some tightness and regularity properties on the sequence 
$\{\nabla\psi^\tau\}$, therefore we apply some results concerning local smoothing due to Vega 
\cite{V} and Constantin, Saut \cite{CS}.
 \begin{prop}[Local smoothing for $\nabla\psi^\tau$]\label{prop:28}
Let $0<T<\infty$ and let $\psi^\tau$ be defined as in the previous section. Then one has
\begin{equation}
\|\nabla\psi^\tau\|_{L^2([0,T];H^{1/2}_{loc}(\real^3))}\leq C(E_0, T, \|\rho_0\|_{L^1}).
\end{equation}
\end{prop}
\begin{dimo}
Using the Strichartz estimates obtained above, we can apply the Theorems \ref{thm:smooth1}, 
\ref{thm:smooth2} about local smoothing. Indeed, by using again the formula \eqref{eq:131} it follows
\begin{align*}
\|\nabla\psi^\tau\|_{L^2([0,T];H^{1/2}_{loc}(\real^3))}\lesssim
&\|\nabla\psi_0\|_{L^2(\real^3)}\\
&+\frac{\tau}{\hbar}\sum_{k=1}^N\|\Lambda^\tau(k\tau-)\|_{L^2(\real^3)}\\
&+\tau\sum_{k=1}^N\|\nabla\psi_{n_k}\|_{L^2(\real^3)}\\
&+\sum_{k=1}^N\norma{\nabla\psi_{n_k}-\nabla\psi^\tau(k\tau-)
+\frac{\tau}{\hbar}(\Lambda_{n_k}-\Lambda^\tau(k\tau-))}_{L^2(\real^3)}\\
&+\|F\|_{L^1([0,T];L^2(\real^3))}.
\end{align*}
The first three terms are clearly estimated by a constant $C(E_0,T)$ depending only on the initial energy and on time. The fourth term is $O(\tau)$. The last term will be estimated using the previous Strichartz estimates. As before, we split $F$ into three parts $F=F_1+F_2+F_3$, then one has
\begin{multline*}
\lplq{|\psi^\tau|^{p-1}\nabla\psi^\tau}{1}{2}{[0,T]}\\\leq T^{\frac{4}{5-p}}
\lplq{|\psi^\tau|^{p-1}}{\frac{4}{p-1}}{\infty}{[0,T]}
\lplq{\nabla\psi^\tau}{\infty}{2}{[0,T]}\\
\leq T^{\frac{4}{5-p}}\lplq{\psi^\tau}{4}{\infty}{[0,T]}^{(p-1)}\lplq{\nabla\psi^\tau}{\infty}{2}{[0,T]},
\end{multline*}
while
\begin{multline*}
\lplq{\nabla V^\tau\psi^\tau}{1}{2}{[0,T]}
\leq T^{\mez}\lplq{\nabla V^\tau}{\frac{2}{1-2\eps}}{\1{\eps}}{[0,T]}
\lplq{\psi^\tau}{\frac{2}{3\eps}}{\frac{2}{1-2\eps}}{[0,T]}
\end{multline*}
and now the remaining calculations are similar to those already done for the Strichartz estimates. Regarding the term 
$V^\tau\nabla\psi^\tau$ we already estimated its $L^1_tL^2_x$ norm.
\end{dimo}
Since $H^{1/2}_{loc}$ is compactly embedded in $L^2_{loc}$, we can apply the Theorem \ref{thm:comp} due to Rakotoson, Temam \cite{RT}.
\begin{prop}\label{prop:29}
The sequence $\{\nabla\psi^\tau\}$ is strongly convergent in 
\newline
$L^2([0, T]; L^2_{loc}(\real^3))$, namely
\begin{equation}\label{eq:221}
\nabla\psi:=s-\lim_{k\to\infty}\nabla\psi^{\tau_k}\qquad\textrm{in}\;L^2([0, T]; L^2_{loc}(\real^3)).
\end{equation}
In particular, one has $\nabla\sqrt{\rho^\tau}\to\nabla\sqrt{\rho}$ and 
$\Lambda^\tau\to\Lambda$ in $L^2([0, T]; L^2_{loc}(\real^3))$.
\end{prop}
\begin{dimo}
The previous Proposition \ref{prop:28} implies that the sequence $\{\nabla\psi^\tau\}_{\tau>0}$, is uniformly bounded in 
$L^2([0,T];H^{1/2}_{loc}(\real^3))$ and then $\nabla\psi^\tau\wlim\nabla\psi$ in that space. Now  
$H^{1/2}_{loc}$ is compactly embedded in $L^2_{loc}$, since $\nabla\psi^\tau(t)\wlim\nabla\psi(t)$ for almost every $t\geq0$ and
\begin{equation*}
\lim_{|E|\to0, E\subset[0,T]}\sup_{\tau>0}\int_E\|\nabla\psi^\tau(t)\|^2_{L^2_{loc}}\de t=0,
\end{equation*}
since $\nabla\psi^\tau\in L^\infty([0,T];L^2(\real^3))$. Hence we can apply the Theorem \ref{thm:comp} by Rakotoson and Temam \cite{RT} and we get \eqref{eq:221}.
\end{dimo}
\begin{prop}
$(\rho, J)$ is a weak solution to the Cauchy problem \eqref{eq:QHD3}, \eqref{eq:QHD_IV3}.
\end{prop}
\begin{dimo}
It follows directly by combining the Theorem \ref{thm:cons} in the section \ref{sect:fract} and the Proposition \ref{prop:29}. As for the collisionless QHD system we should note that the generalized irrotationality condition holds by the definition of the current density and Corollary \ref{cor:irrot}.
\end{dimo}
\subsection*{Acknowledgement}
The authors wish to thank prof. Luigi Ambrosio for some useful comments.

\end{document}